\documentstyle[amssymb,twoside,exti]{article}

\input tcilatex

\newtheorem{hypothesis}[theorem]{Hypothesis}
\def\proof{\noindent{\bf Proof\quad}}

\newcommand{\bP}{{\sf I\hspace*{-.3ex}P}}

\newcommand{\aar}{-\hspace*{-4pt}-\hspace{-5pt}\longrightarrow}

\begin{document}

\author{Donald A.\ Dawson\thanks{%
\thinspace Supported by an NSERC grant and a Max Planck award.}$^{\;\;,}$%
\thanks{%
\thinspace The Fields Institute, 222 College Street, Toronto Ontario, Canada
M5T\thinspace 3J1.\quad E-mail: ddawson@fields.utoronto.ca } \and Klaus
Fleischmann\thanks{%
\thinspace Weierstrass Institute for Applied Analysis and Stochastics
(WIAS), Mohrenstr.\ 39, \quad \quad D--10117 Berlin, Germany.\quad E-mail:
fleischmann@wias-berlin.de, Fax: +49-30-204\thinspace 49\thinspace 75} \and %
Carl Mueller\thanks{%
\thinspace Department of Mathematics, University of Rochester,\quad
Rochester, NY 14620, USA.\quad E-mail: cmlr@troi.cc.rochester.edu}$^{\;\;,}$%
\thanks{%
\thinspace Supported by an NSA grant.}\vspace{20pt}}
\title{Finite time extinction of \vspace{5pt}super-Brownian motions with catalysts%
\vspace{10pt}}
\date{{\scriptsize (WIAS preprint No.\ 431 of 1998)}}
\maketitle

\vfill

\setcounter{page}{0}\thispagestyle{empty}\noindent {\small {\em AMS Subject
Classification}\quad Primary 60J80; Secondary 60J55, 60G57}\bigskip


\vfill 

\noindent {\small Running head:\quad Extinction in catalytic branching}

\vfill 

{\small \ }\hfill {\scriptsize \today\quad exti.tex\quad typeset by \LaTeX}

\vfill\newpage \vspace{-10pt}

\begin{center}
{\small {\bf Abstract}}
\end{center}

\begin{quote}
{\small Consider a catalytic super-Brownian motion $\,X=X^\Gamma $\thinspace
\ with finite variance branching. Here ``catalytic'' means that branching of
the reactant $\,X$\thinspace \ is only possible in the presence of some
catalyst. Our intrinsic example of a catalyst is a stable random measure
\thinspace $\Gamma $\thinspace \ on $\,{\sf R}$\thinspace \ of index $%
\,0<\gamma <1.$\thinspace \ Consequently, here the catalyst is located in a
countable dense subset of $\,{\sf R}.$\thinspace \ Starting with a finite
reactant mass $\,X_0$\thinspace \ supported by a compact set, $\,X$%
\thinspace \ is shown to die in finite time. Our probabilistic argument uses
the idea of good and bad historical paths of reactant ``particles'' during
time periods $\,[T_{n\,},T_{n+1}).$\thinspace \ Good paths have a
significant collision local time with the catalyst, and extinction can be
shown by individual time change according to the collision local time and a
comparison with Feller's branching diffusion. On the other hand, the
remaining bad paths are shown to have a small expected mass at time $%
\,T_{n+1}$\thinspace \ which can be controlled by the hitting probability of
point catalysts and the collision local time spent on them.}
\end{quote}

{\small 
\tableofcontents%
}

\section{Introduction}

Recently a number of papers have dealt with branching in catalytic media.
These are models of chemical or biological reaction-diffusion systems of two
substances or species, respectively. One we call the catalyst, and the other
the reactant. The latter we model by a super-Brownian motion (SBM) with
``critical binary'' branching, and its branching rate is given by the
catalyst.

In this paper we verify {\em finite time extinction}\/ of the reactant for 
{\em three different types of catalysts,}\/ provided the reactant was
started with a finite mass. We begin with explaining the most interesting of
these catalysts.

\subsection{Model 1:\quad stable catalyst $\Gamma $ on ${\sf R}$\label%
{SS.Model1}}

Let $\,X^\Gamma =\{X_t^\Gamma :\,\,t\ge 0\}$\thinspace \ denote a {\em %
continuous super-Brownian motion (SBM)}\/ with branching rate (catalyst)
given by a {\em stable}\/ random measure $\,$%
\begin{equation}
\Gamma \,=\,\sum_{\,i}\alpha _{i\,}\delta _{b_i}
\end{equation}
on the real line $\,{\sf R}$\thinspace \ with index $\,0<\gamma <1.$

At an {\em intuitive level,}\/ this model can be explained as follows. A
huge number of small reactant ``particles'' move independently according to
Brownian motions in ${\sf R}.$ But if such a Brownian particle meets one of
the point catalysts $\,\alpha _{i\,}\delta _{b_i}$ it may branch in a
critical binary way. More precisely, branching is governed by the {\em %
collision local time} 
\begin{equation}
L_{[W,\Gamma ]}({\rm d}s)\,:=\,{\rm d}s\;\sum_{\,i}\alpha _{i\,}\delta
_{b_i}(W_s)
\end{equation}
in the sense of Barlow et al.\ \cite{BarlowEvansPerkins91} of the Brownian
reactant particle with path $\,W$\thinspace \ and the stable random measure $%
\,\Gamma ({\rm d}b)$\thinspace \ describing the catalyst.

This process $\,X^{\Gamma \,}$ was introduced as a Markov process by Dawson
and Fleisch{\Large \-}mann 
\cite[Lemma 2.3.5 and its application in Subsections 2.4\thinspace --\thinspace 2.5]
{DawsonFleischmann91.hfm}. The existence of a continuous version follows
from \cite[Theorem 1\thinspace b)]{DawsonFleischmann97.persist}. The
clumping features of $\,X^{\Gamma \,}$\ had been exhibited in \cite
{DawsonFleischmann91.hfm} by a {\em time-space-mass scaling limit theorem.}%
\/ In \cite{DawsonFleischmannRoelly91.van} the states $\,X_t^{\Gamma \,}$ of 
$\,X^{\Gamma \,}$ had been shown to be {\em absolutely continuous}\/
measures. Finally, in \cite{DawsonLiMueller95}, the so-called {\em compact
support property}\/ has been verified: If the finite initial measure $%
\,X_0^\Gamma $\thinspace \ has compact support, then the range of $%
\,X^{\Gamma \,}$ is compact, too.\smallskip 

Starting with a finite measure $\,X_{0\,}^\Gamma ,$\thinspace \ and given $%
\,\Gamma ,$\thinspace \ the total mass process $\,t\mapsto X_t^\Gamma ({\sf R%
})$\thinspace \ is a continuous martingale, hence has a.s.\ a limit as $%
\,t\uparrow \infty $\thinspace \ (\cite[Proposition 3]
{DawsonFleischmann97.persist}). The {\em main purpose}\/ of the present
paper is to show that if $\,X_{0\,}^\Gamma $ is of compact support, the
process $\,X^\Gamma $\thinspace \ {\em dies in finite time}\/ (Theorem \ref
{T.finite.time.exti} on page \pageref{T.finite.time.exti}), just as in the
constant medium case (the formal $\,\gamma =1$\thinspace \ boundary
case).\smallskip 

To illustrate the problems we encounter in the proof, we consider the
following. Given the catalyst $\,\Gamma $\thinspace \ and starting $%
\,X^{\Gamma \,}$ with a unit mass concentrated at $\,a,$\thinspace \ that is 
$\,X_0^\Gamma =\delta _{a\,},$\thinspace \ the probability of extinction of $%
\,X^{\Gamma \,}$ at time $\,t$\thinspace \ is given by 
\begin{equation}
P_{0,\delta _a}^\Gamma (X_t^\Gamma =0)\;=\;\exp \left[ -v_\infty \big(0,a\,%
\big| \,t,\Gamma \big)\right]   \label{exti.prob}
\end{equation}
where for $\,\theta ,t,\Gamma $\thinspace \ fixed, $v_\theta =\,v_\theta
(\cdot ,\cdot \,|\,t,\Gamma )=\big\{v_\theta (s,a\,|\,t,\Gamma ):\,(s,a)\in
[0,t]\times {\sf R}\big\}$ solves (formally) the following {\em %
reaction-diffusion equation in the stable catalytic medium} $\,\Gamma :$ 
{\em \ } 
\begin{equation}
-\,\frac \partial {\partial s}\,v_\theta \,=\,\frac 12\,\Delta v_\theta
-\Gamma \,v_\theta ^2,\vspace{4pt}\qquad {\rm \,}v_\theta (s,a\,|\,t,\Gamma )%
\big| _{s=t}\,\equiv \,\theta \ge 0,  \label{equ}
\end{equation}
and $\,v_\infty :=\lim\nolimits_{\theta \uparrow \infty \,}v_{\theta \,}.$%
\thinspace \ Then, by Borel-Cantelli, it would suffice to show the following
extinction property of solutions to (\ref{equ}): 
\begin{equation}
\lim_{t\uparrow \infty }\,\lim_{\theta \uparrow \infty }\,v_\theta
(0,a\,|\,t,\Gamma )\,=\,0.  \label{solution.vanishes}
\end{equation}
But we do not know how to attack this problem analytically. Recall that the
coefficient $\,\Gamma $\thinspace \ of the reaction term (reaction rate) in (%
\ref{equ}) is the generalized derivative of a (random) measure supported by
a countable dense set in $\,{\sf R}$,\thinspace \ hence is highly singular.

Instead, to prove finite time extinction will use some {\em probabilistic}\/
arguments concerning the stochastic process $\,X^\Gamma ,\,\,\;$inspired by
Fleischmann and Mueller \cite{FleischmannMueller97.infcatmass}.

At the same time, via the log-Laplace connection of $\,X^{\Gamma \,}$ to the
partial differential equation (\ref{equ}), our approach can be regarded as a
probabilistic contribution to the study of asymptotic properties [such as (%
\ref{solution.vanishes})] of solutions to the reaction-diffusion equation (%
\ref{equ}) in the (random) heterogeneous singular medium $\,\Gamma .$%
\thinspace \ 

Equations as (\ref{equ}) have attracted some attention and are relevant in
particular from an {\em applied}\/ point of view; see e.g.\ Ortoleva and
Ross \cite{OrtolevaRoss72}, Pagliaro and Taylor \cite{PagliaroTaylor88}. For
reaction-diffusion equations in heterogenous media with different species
and where reaction may be concentrated on bounded interfaces, see Glitzky et
al.\ \cite{GlitzkyGroegerHuenlich96}. Note that reaction-diffusion equations
arise in many branches of technology, e.g.\ in microelectronics.\smallskip

The {\em main ideas}\/ of our approach are as follows. First of all, since
we start with an initial measure $\,X_0$\thinspace \ of compact support, and 
$\,X^\Gamma \,\,\;$has the compact support property (\cite{DawsonLiMueller95}%
), we may ``essentially'' restrict our attention to a finite (space)
interval $\,{\sf K}\subset {\sf R}.$\thinspace \ Hence, by a coupling
technique, the catalyst may be extended periodically outside $\,{\sf K}.$%
\thinspace \ Next, the probability of extinction can be estimated below by
using a smaller branching rate. Therefore, we remove all atoms $\,\alpha
_{i\,}\delta _{b_i}$\thinspace \ of the catalyst with large ``action
weights'' $\,\alpha _{i\,}.$\thinspace \ Moreover, the action weights $%
\,\alpha _i$\thinspace \ belonging to $\,[2^{-n},2^{-n+1})$\thinspace \ are
replaced by $\,2^{-n},$\thinspace \ so that the corresponding atoms form
Poisson point processes in $\,{\sf K}$\thinspace \ of intensity $%
\,c\,2^{\gamma n}$\thinspace \ (with $\,c$\thinspace \ an appropriate
constant). Finally, ``large'' distances between neighboring points of this
Poisson point process are exceptional. Therefore we may restrict to the
situation where the empty intervals are at most of a size $\,\Delta _n$%
\thinspace \ (to be specified later). Altogether, we then want to verify
finite time extinction of $\,X^{\underline{\Gamma }}$\thinspace \ where the
catalyst $\,\underline{\Gamma }$\thinspace \ is of the form $\,\sum_{n\ge
0\,}2^{-n}\underline{\pi }_{\,n}\,\;$\vspace{1pt}where $\,\underline{\pi }%
_{\,n}$\thinspace \ is a periodic point measure with gaps between
neighboring catalysts bounded by $\,\Delta _{n\,}.$

The {\em central idea}\/ is to look for a sequence of times $\,T_1<T_2<\cdot
\cdot \cdot $\thinspace \ with {\em finite}\/ accumulation point $\,T_\infty 
$\thinspace \ with the following property. At time $\,T_{n\,},$\thinspace \
we distinguish between \vspace{1pt}{\em ``good and bad''}\/ historical paths
of Brownian reactant particles, starting from the state $\,X_{T_n}^{%
\underline{\Gamma }}$\thinspace \ at time $\,T_{n\,},$\thinspace \ as we now
explain.

The {\em good}\/ paths are those which have a ``significant'' collision
local time with $\,2^{-n}\underline{\pi }_{\,n}$\thinspace \ on the time
interval $\,[T_{n\,},T_{n+1})$\thinspace \ (so we take into account only
that part $\,2^{-n}\underline{\pi }_{\,n}$\thinspace \ of $\,\underline{%
\Gamma }).$\thinspace \ Consider the total mass of the good paths. For the
continuous SBM with a {\em uniform}\/ branching rate,\ the total mass
process would have the same distribution as the standard {\em Feller
branching diffusion}\/ which satisfies the one-dimensional stochastic
equation 
\begin{equation}
{\rm d}Z_r\,=\,\sqrt{2\,Z_r}\,{\rm d}B_{r\,},\quad Z_0\ge 0,
\label{Feller's.diffusion}
\end{equation}
(with $\,B$\thinspace \ a standard Brownian motion). It is well-known that
this diffusion is absorbed at $\,0$\thinspace \ in finite time. In our
catalytic case, the total mass of the good paths can in law be {\em compared}%
\/ with Feller's branching diffusion. But now its time scale during $%
\,[T_{n\,},T_{n+1})$\thinspace \ is, roughly speaking, individually given by
the collision local times of the good paths with the catalytic medium $%
\,2^{-n}\underline{\pi }_{\,n\,}.$\thinspace \ Since these collision local
times are ``significant'' on the good paths, it follows that the total mass
of the good paths dies out by time\ $\,T_{n+1}$\thinspace \ with high
probability.

The remaining {\em bad}\/ paths may {\em not}\/ die out by time $\,T_{n+1\,},
$\thinspace \ but we can estimate the probability that this mass is larger
than a certain size at time $\,T_{n+1\,},$\thinspace \ by using Markov's
inequality and the simple but powerful expectation formula for (historical)
superprocesses. Then we need to derive some estimates concerning hitting
probabilities of a neighboring point catalyst from $\,2^{-n}\underline{\pi }%
_{\,n\,},$\thinspace \ and the Brownian (collision) local time spent on it.

\subsection{Model 2:\quad i.i.d.\ uniform catalysts on the lattice ${\sf Z}%
^d $\label{SS.Model2}}

In the other two models we discuss, the basic ideas of distinguishing
between good and bad historical paths, and how to handle them, are the same.
So here we only introduce the models, and indicate how to classify the paths.

For the second model, we replace the phase space $\,{\sf R}$\thinspace \ by
the lattice $\,{\sf Z}^d,$\thinspace \ and Brownian motion by a continuous
time {\em simple random walk\/} in $\,{\sf Z}^d.\,\,\;$The catalysts $%
\,\varrho =\{\varrho _b:\,b\in {\sf Z}^d\}$\thinspace \ are {\em i.i.d.}\/\
random variables, {\em uniform}\/ in the interval $\,(0,1).$\thinspace \ So
here the catalysts are present everywhere but again their action weights
fluctuate randomly. By the discreteness of $\,{\sf Z}^d,$ and since masses
can be arbitrarily small in superprocesses, one does not expect that the
compact support property holds. Therefore, as opposed to Model 1, the
super-random walk $\,X^\varrho $\thinspace \ with catalyst $\,\varrho $%
\thinspace \ may be influenced by large regions, where the catalysts are
small. However, calling paths {\em bad}\/ which reach such a region, these
paths should have a small expected mass, and we will be able to show the
finite time extinction property for $\,X^\varrho $\thinspace \ along the
lines indicated.

\subsection{Model 3:\quad a deterministic ``parabolic'' catalyst $\chi _q$ 
\label{SS.Model3}}

For the moment, consider the continuous SBM with phase space ${\sf R}$ and
uniform branching rate, except on $\,(-1,1).$\thinspace \ More precisely, we
consider the branching rate $\,\chi ={\bf 1}_{{\sf R}\backslash (-1,1)\,}.$%
\thinspace \ As we will see in the next Subsection, if $\,X_0\big((-1,1)\big)%
>0,$\thinspace \ then this superprocess does not die in finite time.

Motivated by this, for a fixed constant $\,q>0,$\thinspace \ we consider the
truncated {\em ``parabolic''}\/ branching rate 
\begin{equation}
\chi _q(b)\,:=\,|b|^q\wedge 1,\qquad b\in {\sf R},  \label{parabolic.rate}
\end{equation}
(see the figure). We show that starting with a finite initial mass, the SBM $%
\,X^\chi $\thinspace \ with parabolic catalyst $\,\chi _{q\,}$\ {\em dies in
finite time,}\/ just as in the constant branching rate case, despite the
``depression'' of branching rate close to the origin, even if $\,q$%
\thinspace \ is very large. Here the good historical paths are those which
do not spent too much time near $\,0,\,\,\;$where the catalytic mass is
small. 
\[
\fbox{$\vspace{3pt}
\begin{array}{ccc}
\FRAME{itbpF}{3.0006cm}{1.0017cm}{0cm}{}{}{Plot }{\special{language
"Scientific Word";type "MAPLEPLOT";width 3.0006cm;height 1.0017cm;depth
0cm;display "USEDEF";plot_snapshots TRUE;function \TEXUX{$\left| b\right|
^{1/3}\wedge 1$};linecolor "black";linestyle 1;linethickness 1;pointstyle
"point";xmin "-4";xmax "4";xviewmin "-4";xviewmax "4";yviewmin
"0.1364";yviewmax "1.017";rangeset"X";recompute TRUE;phi 45;theta
45;plottype 4;numpoints 50;axesstyle "none";xis \TEXUX{b};var1name
\TEXUX{$x$};valid_file "T";tempfilename
'EYASU401.wmf';tempfile-properties "XP";}}\vspace{2pt} & 
\FRAME{ipF}{3.0006cm}{1.0017cm}{0cm}{}{}{Plot }{\special{language
"Scientific Word";type "MAPLEPLOT";width 3.0006cm;height 1.0017cm;depth
0cm;display "USEDEF";plot_snapshots TRUE;function \TEXUX{$\left| b\right|
\wedge 1$};linecolor "black";linestyle 1;linethickness 1;pointstyle
"point";xmin "-4";xmax "4";xviewmin "-4";xviewmax "4";yviewmin
"-0.01632";yviewmax "1.02";rangeset"X";recompute TRUE;phi 45;theta
45;plottype 4;numpoints 50;axesstyle "none";xis \TEXUX{b};var1name
\TEXUX{$x$};valid_file "T";tempfilename
'EYASVW02.wmf';tempfile-properties "XP";}}\vspace{2pt} & 
\FRAME{ipF}{3.0006cm}{1.0017cm}{0cm}{}{}{Plot }{\special{language
"Scientific Word";type "MAPLEPLOT";width 3.0006cm;height 1.0017cm;depth
0cm;display "USEDEF";plot_snapshots TRUE;function \TEXUX{$\left| b\right|
^5\wedge 1$};linecolor "black";linestyle 1;linethickness 1;pointstyle
"point";xmin "-4";xmax "4";xviewmin "-4";xviewmax "4";yviewmin
"-0.02";yviewmax "1.02";rangeset"X";recompute TRUE;phi 45;theta 45;plottype
4;numpoints 70;axesstyle "none";xis \TEXUX{b};var1name
\TEXUX{$x$};valid_file "T";tempfilename
'EYASWV03.wmf';tempfile-properties "XP";}}\vspace{2pt} \\ 
\chi _{1/3}(b)=\left| b\right| ^{1/3}\wedge 1 & \chi _1(b)=\left| b\right|
\wedge 1\vspace{12pt} & \chi _5(b)=\left| b\right| ^5\wedge 1 \\ 
& \hspace*{-15pt}{\rm Variants\;of\;the\;``parabolic"\;catalyst}%
\hspace*{-15pt} & 
\end{array}
$}
\]

\subsection{Non-extinction in finite time\label{SS.nonextinction}}

If we change Model 1 so that the catalysts are not\/ dense, then the mass
fails to die out in finite time. In fact, if $\,I\ne \emptyset $\thinspace \
is an open interval without catalysts, then a corresponding catalytic SBM $%
\,X$\thinspace \ is bounded below by the heat flow in $\,I$\thinspace \ with
absorption at the boundary $\,\partial I,$\thinspace \ starting with $%
\,X_0(\cdot \,\cap I).$\thinspace \ If now $\,X_0(I)>0,$\thinspace \ then
the $L^1$--norm of that heat solution decays according to $\;\left\langle
X_{0\,},\varphi _\lambda \right\rangle \,{\rm e}^{-\lambda t}$\thinspace \
with $\,\lambda >0$\thinspace \ the first eigenvalue of $\,\frac 12\,\Delta $%
\thinspace \ on $\,I,$\thinspace \ and $\,\varphi _\lambda $\thinspace \ is
the corresponding eigenfunction, hence is (strictly) positive at any time $%
\,t,$\thinspace \ that is $\,X_t(I)>0$\thinspace \ for all $\,t.$

Note that catalytic SBMs with a gap cover the single point catalytic SBM,
where survival for all finite times was known from Fleischmann and
Le\thinspace Gall \cite{FleischmannLeGall95}.\smallskip

It would be interesting to establish conditions on the catalytic medium
which are necessary {\em and}\/ sufficient for extinction in finite time.
Unfortunately, our methods seem to be too crude for this.

\begin{remark}[decomposition of initial measures]
\label{R.decomposition}{\em Suppose a decomposition }$\,\mu =\sum_i\mu _i$%
{\em \thinspace \ of the initial measure is given. If we can show finite
time extinction for each initial measure }$\,X_0=\mu _i$\thinspace \ {\em %
then the branching property implies finite time extinction for }$\,X_0=\mu .$%
\hfill $\Diamond $
\end{remark}

\subsection{Outline\label{SS.outline}}

To give a precise meaning to the above ideas, some technical problems have
to be overcome. For instance, \vspace{1pt}to have access to reactant
particle paths, we will work with the {\em historical}\/ catalytic SBM $%
\widetilde{X}^\Gamma $ instead of $\,X^\Gamma .$\thinspace \ Or, since we
want to use time scales of individual reactant particles, we will exploit
Dynkin's \cite{Dynkin91.a} framework of ``stopped'' historical
superprocesses.\smallskip 

The {\em outline}\/ of the paper is as follows. In the next section we
recall the model of {\em continuous SBM}\/ $X$ in $\,{\sf R}^d$\thinspace \
with branching rate functional $\,K$\thinspace \ as provided in Dawson and
Fleischmann \cite{DawsonFleischmann97.persist} (this goes back to Dynkin 
\cite{Dynkin91.a}). Then $\,K$\thinspace \ is specialized to be the Brownian
collision local time $\,K=L_{[W,\psi ]}$\thinspace \ of a (deterministic)
locally finite measure $\,\psi $\thinspace \ (catalyst) on $\,{\sf R},$%
\thinspace \ also taken from \cite{DawsonFleischmann97.persist}. Further
specialized to $\,\psi =\Gamma ,$\thinspace \ our main result, Theorem \ref
{T.finite.time.exti} \vspace{1pt}at p.\pageref{T.finite.time.exti}, can be
formulated.\vspace{1pt}

In Section \ref{S.general.criterion} we first recall the {\em historical}\/
SBM $\,\widetilde{X}$\thinspace \ in $\,{\sf R}^d$\thinspace \ with
branching rate functional $\,K.$\thinspace \ For this model, we give an {\em %
abstract sufficiency criterion}\/ (Theorem \ref{T.abstract.criterion} at p.%
\pageref{T.abstract.criterion}) for finite time extinction based on the idea
of good and bad paths. For the extinction of good paths, a {\em comparison}%
\/ with Feller's branching diffusion is provided (Proposition \ref
{P.comparison.Feller} at p.\pageref{P.comparison.Feller}), as a refinement
of an argument in \cite{FleischmannMueller97.infcatmass}. This makes use of
Dynkin's concept of (individually) {\em ``stopped''}\/ historical
superprocesses.

Section \ref{S.application} is devoted to two one-dimensional {\em %
applications}\/ of the abstract criterion:\vspace{-7pt}

\begin{description}
\item[(i) ]  the parabolic catalyst $\,\chi _q$\thinspace \ of Model 3, and%
\vspace{-5pt}

\item[(ii)]  \hspace*{2pt}a (deterministic) point catalyst $\,\vspace{2pt}%
\underline{\Gamma }=\sum_{n\ge N\,}2^{-n}\underline{\pi }_{\,n}$\thinspace \
with dense locations and gaps between neighboring catalysts in $\,\underline{%
\pi }_{\,n}$\thinspace \ bounded by some $\,\Delta _{n\,}$.\vspace{2pt}%
\vspace{-8pt}
\end{description}

In Section \ref{S.proof.main} we prove our main theorem, the finite time
extinction for the SBM $\,X^\Gamma $\thinspace \ with a stable catalytic
rate $\,\Gamma $\thinspace \ (Model 1). In fact, by a coupling and
comparison argument, we reduce the problem to the case (ii) above.

Finally, in Section \ref{SS.lattice}, finite time extinction for the
super-random walk $\,X^\varrho $\thinspace \ with i.i.d.\ uniform catalysts
is derived.

\section{Stable catalysts -- main result\label{S.main}}

Here we carefully introduce the continuous SBM $\,X$\thinspace \ in $\,{\sf R%
}^d$\thinspace \ with a sufficiently nice branching rate functional $\,K.$%
\thinspace \ After specializations to Model 1, we will formulate our main
result, Theorem \ref{T.finite.time.exti} at p.\pageref{T.finite.time.exti}.

\subsection{Preliminaries: some spaces\label{SS.preliminaries}}

Measurability is always meant with respect to the related Borel fields. The
lower index + refers to the subset of all non-negative members of a set.

Let $\,{\cal B}[E_{1\,},E_2]$\thinspace \ denote the set of all {\em %
measurable}\/ mappings $\,f:E_1\rightarrow E_2$\thinspace \ where $%
\,E_{1\,},E_2$\thinspace \ are topological spaces. Write $\,{\cal B}[E_1]$%
\thinspace \ instead of $\,{\cal B}[E_{1\,},E_2]$\thinspace \ if $\,E_2={\sf %
R},$\thinspace \ the real line, and only $\,{\cal B}$\thinspace \ if
additionally $\,E_1={\sf R}^d,$ $\,\,d\ge 1.$\thinspace \ 

If we restrict our consideration to {\em continuous\/} functions $\,f,$%
\thinspace \ the letter $\,{\cal B}$\thinspace \ is replaced by $\,{\cal C}$%
\thinspace \ in the respective cases. If we restrict to {\em bounded}\/
functions, we write $\,b{\cal B}$\thinspace \ and $\,b{\cal C},$\thinspace \
etc.

Fix a dimension $\,d\ge 1,$\thinspace \ and a constant $\,p>d,$\thinspace \
and introduce the {\em reference function} 
\begin{equation}
\phi _p(b)\,:=(1+|b|^2)^{-p/2},\qquad b\in {\sf R}^d,
\label{reference.function}
\end{equation}
of $\,p$--potential decay at infinity. Denote by $\,{\cal B}^{\,p}$%
\thinspace \ the set of all those $\,\varphi \in {\cal B}$\thinspace \ such
that $\,|\varphi |\le c_\varphi \phi _p$\thinspace \ for some constant $%
c_{\varphi \,}.$\vspace{1pt}

Write $\,\left\langle \mu ,f\right\rangle $\thinspace \ for the integral $%
\,\int \!\mu ({\rm d}b)\,f(b).$\thinspace \ Let $\,{\cal M}_p={\cal M}_p[%
{\sf R}^d]\,\;$denote the set of all (non-negative) measures $\,\mu $%
\thinspace \ defined on $\,{\sf R}^d$\thinspace \ satisfying $\,\left\langle
\mu ,\phi _p\right\rangle <\infty .$\thinspace \ We endow this set $\,{\cal M%
}_p$\thinspace \ of $\,p$--tempered measures with the weakest topology such
that all the maps $\,\mu \mapsto \left\langle \mu ,\varphi \right\rangle $%
\thinspace \ are continuous, where $\,\varphi \ge 0$\thinspace \ is
continuous and of compact support, or $\,\varphi =\phi _{p\,}.$\thinspace \
The set of all finite measures on a Polish space $\,E$\thinspace \ is
denoted by $\,{\cal M}_{{\rm f}}[E]$\thinspace \ and equipped with the
topology of weak convergence. Write $\,\Vert {}\mu \Vert \,\;$for the total
mass $\,\mu (E)=\left\langle \mu ,1\right\rangle \,\;$of $\,\mu \in {\cal M}%
_{{\rm f}}[E].$\vspace{1pt}

Set $\,{\cal M}_{{\rm f}}={\cal M}_{{\rm f}}\left[ {\sf R}^d\right] ,$%
\thinspace \ and denote by $\,{\cal M}\left[ {\sf R}_{+}\times E\right] $%
\thinspace \ the set of all measures $\,\eta $\thinspace \ defined on $\,%
{\sf R}_{+}\times E$\thinspace \ such that $\,\eta \left(
[0,t]_{\!_{\!_{\,}}}\times E\right) <\infty $\thinspace \ for all $\,t>0.$

With $\,c$\thinspace \ we always denote a positive constant which may be
different at various places. On the other hand, constants $\,c_i$\thinspace
\ are fixed within each subsection.

\subsection{Branching rate functional $K$ and BCLT $L_{[W,\psi ]}$\label%
{SS.BCLT}}

Let $\,W=[W,\,\Pi _{s,a\,},\,s\ge 0,\,\;a\in {\sf R}^d]$\thinspace \ denote
the {\em standard Brownian motion} in $\,{\sf R}^d,$\thinspace \ on
canonical path space $\,{\cal C}\left[ {\sf R}_{+\,},{\sf R}^d\right] $%
\thinspace \ of continuous functions $\,w.$

\begin{remark}[inhomogeneous setting]
\label{R.inhomogeneous}{\em Although Brownian motion is time-homogeneous, we
use this inhomogeneous setting, and we read }$\,\Pi _{s,a\,}\varphi (W_t)$%
\thinspace \ {\em as }$\,0$\thinspace \ {\em if }$\,s>t.$\thinspace \ {\em %
This formalism looks artificial since changing the paths before time }$\,s$%
\thinspace \ {\em does not change the laws }$\,\Pi _{s,a\,}.$\thinspace \ 
{\em The advantage becomes clear when we work with historical SBM. Note that
the measure }$\,\Pi _{s,a}$\thinspace \ {\em is concentrated on the set of
paths }$\,\left\{ w\in {\cal C}[{\sf R}_{+\,},{\sf R}^d]:_{\!_{\!_{\,}}}%
\,w_s=a\right\} .$\hfill $\Diamond $
\end{remark}

Write $\,{\rm p}$\thinspace \ for the continuous {\em transition density}\/
function of $\,W$, 
\begin{equation}
{\rm p}_t(a,b)\,=\,{\rm p}_t(b-a)\,=\,(2\pi t)^{-d/2}\,\exp \left[ -\frac{%
(b-a)^2}{2t}\right] ,\quad \,\;t>0,\,\;a,b\in {\sf R}^d,
\label{Brownian.transition.densities}
\end{equation}
and 
\begin{equation}
\Pi _\eta \,:=\,\int \!\!\eta ({\rm d}s,{\rm d}a)\,\Pi _{s,a\,},\qquad \eta
\in {\cal M}\left[ {\sf R}_{+}\times {\sf R}^d\right] ,
\end{equation}
for the ``law'' of $\,W$ starting at time $\,s$\thinspace \ in a point $\,a$%
\thinspace \ where $\,(s,a)$\thinspace \ is ``distributed'' according to the
measure $\,\eta .$\thinspace \ Put 
\begin{equation}
\Pi _{s,\mu }\,:=\,\Pi _{\delta _s\times \mu \,},\qquad s\ge 0,\quad \mu \in 
{\cal M}_{{\rm f\,}}.
\end{equation}

For convenience, we introduce the following definition.

\begin{definition}[branching rate functional]
\label{D.brf}{\em A continuous additive functional }$\,K=K_{[W]}$\thinspace
\ {\em of Brownian motion }$\,W$\thinspace \ {\em is called a }branching
rate functional\/{\em \ in }$\,{\bf K}^\nu ,$\thinspace \ {\em for some }$%
\,\nu >0,$\thinspace \ {\em if the following two conditions hold:}\vspace{%
-6pt}

\begin{description}
\item[(a) ]  {\em It is} locally admissible,\footnote{$^{)}\,$For
non-admissible functionals, we refer to \cite{DawsonFleischmannLeduc97.cont}.%
}$^{)}$ {\em i.e.\ } 
\[
\sup_{a\in {\sf R}^d\,}\Pi _{s,a}\int_s^t\!K({\rm d}r)\,\phi _p(W_r)\,%
\stackunder{s,t\rightarrow r_0}{%
\aar%
}\,0,\qquad r_0\ge 0.\vspace{-6pt}
\]

\item[(b) ]  {\em For each }$\,N$\thinspace \ {\em there is a constant }$%
\,c_N>0$\thinspace \ {\em such that} 
\[
\Pi _{s,a}\int_s^t\!K({\rm d}r)\,\phi _p^2(W_r)\,\le \,c_N\,|t-s|^\xi \,\phi
_p(a),
\]
$0\le s\le t\le N,\,\;a\in {\sf R}^d.$\hfill $\Diamond $
\end{description}
\end{definition}

To come to our {\em main example}\/ of a branching rate functional, consider
for the moment $\,d=1$\thinspace \ and fix $\,\psi \in {\cal M}_{p\,}.$%
\thinspace \ Intuitively, 
\begin{equation}
L_{[W,\psi ]}({\rm d}r)\;:=\;\,{\rm d}r\,\int \!\!\psi ({\rm d}b)\,\delta
_b(W_r)  \label{BCLT}
\end{equation}
is the {\em Brownian collision local time (BCLT)}\/ of $\,\psi .\,\;$From 
\cite[Corollary 2, p.257]{DawsonFleischmann97.persist} we immediately get
the following statement.\thinspace \footnote{$^{)}$\thinspace For the more
general case if $\varrho \,$ is also time-dependent or, in particular a path
of ordinary SBM, we refer to \cite{FleischmannKlenke98.smooth} and
references therein.}$^{)}$

\begin{lemma}[Brownian collision local time of $\,\varrho $]
\label{L.BCLT.Gamma}Fix\/ $\,d=1$\thinspace \ and\/ $\,\psi \in {\cal M}%
_{p\,}.$\thinspace \ The Brownian collision local time\/ $\,L_{[W,\psi ]}$%
\thinspace \ of\/ $\,\psi $\thinspace \ makes sense non-trivially as a
continuous additive functional of \vspace{1pt}(one-dimensional) Brownian
motion $\,W,$\thinspace \ and it is a branching rate functional in\/ $\,{\bf %
K}^\nu $\thinspace \ with\/ $\,\xi =\frac 12\,.$
\end{lemma}

\subsection{SBM $X$ with branching rate functional $K$\label%
{SS.catalytic.SBM}}

A slight modification of Proposition 12 (p.230) and Theorem 1 (p.234) in 
\cite{DawsonFleischmann97.persist} gives the following lemma.

\begin{lemma}[continuous SBM with branching rate functional $\QTR{bf}{K}$]
\label{L.cont.SBM.stable}\hfill Fix a\newline
dimension\/ $\,d\ge 1,$\thinspace \ and\/ $\,K\in {\bf K}^\nu \,\;$for
some\/ $\,\nu >0.$ \vspace{-5pt}

\begin{description}
\item[(a) (existence) ]  There exists a continuous\/ ${\cal M}_{{\rm f}}$%
--valued (time-inhomogeneous) Markov process\/ $\,X=\big[X,\,P_{s,\mu
\,},\,s\ge 0,\,\mu \in {\cal M}_{{\rm f}}\big]$\thinspace \ with Laplace
functional 
\begin{equation}
P_{s,\mu }\exp \left\langle X_{t_{\,_{}}},-\varphi \right\rangle \,=\,\,\exp
\left\langle \mu \,,_{\!_{\!_{\,}}}-v(s,\cdot \,|\,t)\right\rangle ,
\label{Laplace.X.Gamma}
\end{equation}
$0\le s\le t,\,\;\mu \in {\cal M}_{{\rm f}},\,\;\varphi \in b{\cal B}%
_{+\,},\,\;$where for\/ $\,t,\varphi $\thinspace \ fixed,\/ $\,v=v(\cdot
\,,\cdot \,|\,t)\ge 0$\thinspace \ is uniquely determined by the log-Laplace
equation 
\begin{equation}
v(s,a)\,=\,\Pi _{s,a}\Big[\varphi (W_t)-\int_s^t\!K({\rm d}r)\,v^2(r,W_r)%
\Big] \,,  \label{log.Laplace.X.Gamma}
\end{equation}
$0\le s\le t,\,\;a\in {\sf R}^d.$\vspace{-3pt}

\item[(b) (modification) ]  To each\/ $\,\eta \in {\cal M}\left[ {\sf R}%
_{+}\times {\sf R}^d\right] ,$\thinspace \ there is an\/ ${\cal M}_{{\rm f}}$%
--valued Mar{\Large \-}kov process $\,\left[ X,P_\eta \right] $\thinspace \
such that 
\begin{equation}
P_\eta \exp \left\langle X_{t_{\,_{}}},-\varphi \right\rangle \,=\,\,\exp
\left\langle \eta \,,_{\!_{\!_{\,}}}-v(\cdot ,\cdot \,|\,t)\right\rangle
,\qquad t\ge 0,
\end{equation}
with\/ $\,v(s,\cdot \,|\,t)$\thinspace \ from\/ {\em (a) }if\/ $\,0\le s\le
t,$\thinspace \ and\/ $\,v(s,\cdot \,|\,t)=0$\thinspace \ otherwise.\vspace{%
-4pt}

\item[(c) (moments) ]  $\left[ X,P_{s,\mu }\right] $\thinspace \ \vspace{1pt}%
has finite moments of all orders. In particular, for\/ $\,\eta \in {\cal M}%
\left[ {\sf R}_{+}\times {\sf R}^d\right] $\thinspace \ and\/ $\,\varphi
_1,\varphi _2\in b{\cal B}_{+\,},$\thinspace \ as well as\/ $%
\,t_1,t_{2\,}\ge 0$%
\begin{equation}
P_\eta \left\langle X_{t_1\,},\varphi _1\right\rangle \,=\,\Pi _{\eta
\,}\varphi _1(W_{t_1}),
\end{equation}
\begin{equation}
\left. 
\begin{array}{l}
C{\rm ov}_\eta \big[\!\left\langle X_{t_1\,},\varphi _1\right\rangle
\!,\!\left\langle X_{t_2\,},\varphi _2\right\rangle \!\big]\vspace{5pt} \\ 
=\,2\,\Pi _\eta \dint \!\!K({\rm d}r)\,\left[ \Pi _{r,W_{r_{}}}\varphi
_1\left( W_{t_1}\right) \right] \left[ \Pi _{r,W_{r_{}}}\varphi _2\left(
W_{t_2}\right) \right] .
\end{array}
\;\right\} 
\end{equation}
\end{description}
\end{lemma}

This superprocess $\,X$\thinspace \ is said to be the {\em continuous SBM
with branching rate functional }$\,K.$\thinspace \ 
Note that the lemma in particular applies in the case of a BCLT $%
\,K=L_{[W,\psi ]}$\thinspace \ \vspace{2pt}according to Lemma \ref
{L.BCLT.Gamma}, resulting in a time-homogeneous Markov process.

\subsection{Main result: finite time extinction of $X^\Gamma $\label%
{SS.main.result}}

Let $\,d=1.$\thinspace \ Fix a constant $\,0<\gamma <1,$\thinspace \ and a
(not necessarily normalized) Lebesgue measure $\,\ell $\thinspace \ on ${\sf %
R}.$ The {\em stable catalyst}\/ $\,\left( \Gamma ,%
\bP%
\right) $\thinspace \ is by definition the stable random measure on ${\sf R}$
with Laplace functional 
\begin{equation}
\bP%
\exp \left\langle \Gamma ,-\varphi \right\rangle \;=\;\exp \left[ -\int
\!\!\ell ({\rm d}b)\,\varphi ^\gamma (b)\right] ,\qquad \varphi \in {\cal B}%
_{+\,}.  \label{Laplace.Gamma}
\end{equation}
Recall that $\,\Gamma $\thinspace \ has independent increments, and that it
allows a {\em representation}\/ $\,$%
\begin{equation}
\Gamma \,=\,\sum_{\,i}\alpha _{i\,}\delta _{b_i}  \label{Gamma.rep}
\end{equation}
where the set of locations $\,b_i$\thinspace \ is {\em dense}\/ in $\,{\sf R}%
,$\thinspace \ with $\,%
\bP%
$--probability one.

We now additionally require $\,p>\tfrac 1\gamma $\thinspace \ [for the
exponent $\,p$\thinspace \ of potential decay occurring in the reference
function (\ref{reference.function})]. Then by (\ref{Laplace.Gamma}), the
realizations of the catalyst $\,\Gamma $\thinspace \ belong $%
\bP%
$--almost surely to $\,{\cal M}_{p\,}.$\thinspace \ Hence we may apply the
constructions of the previous two subsections to introduce the {\em %
continuous SBM}\/ $X^\Gamma =\big[X^\Gamma ,\,P_{s,\mu \,}^\Gamma ,\,s\ge
0,\,\mu \in {\cal M}_{{\rm f}}\big]$ {\em with stable catalyst\/ }$\,\Gamma .
$\thinspace \ More precisely, we use the so-called {\em quenched}\/
approach: First a realization $\,\Gamma $\thinspace \ of the catalytic
medium is selected according to $\,%
\bP%
,$\thinspace \ and then, given $\Gamma ,$ the continuous time-homogeneous
Markov process $\,X^\Gamma $\thinspace \ evolves, governed by the BCLT $%
\,L_{[W,\Gamma ]\,}$.\vspace{1pt}

Note that by a formal differentiation of the log-Laplace equation (\ref
{log.Laplace.X.Gamma}) with $\,K=L_{[W,\Gamma ]}$\thinspace \ of (\ref{BCLT}%
), using the semigroup of $\,W,$\thinspace \ and replacing $\,\varphi $%
\thinspace \ by the constant function $\,\theta ,$\thinspace \ we get back
the reaction-diffusion equation (\ref{equ}) in the catalytic medium $%
\,\Gamma .$

Now we are in a position to formulate our {\em main result.}\/ Recall that $%
\,d=1.$

\begin{theorem}[finite time extinction of $X^\Gamma $]
\label{T.finite.time.exti}Fix\/ $\,\mu \in {\cal M}_{{\rm f}}$\thinspace \
with compact support. For\/ $\,%
\bP%
$--almost all\/ $\,\Gamma $\thinspace \ the following holds: 
\begin{equation}
P_{0,\mu }^\Gamma \Big(X_t^\Gamma =0\,\;{\rm for\;some\;\,}t\Big)\,=\,1.%
\vspace{4pt}  \label{finite.time.exti}
\end{equation}
\end{theorem}

The proof of this theorem needs some preparation and is postponed until
Section \ref{S.proof.main}.

We mention that it is an {\em open problem}\/ whether finite time extinction
holds also in some $\,\gamma =0$\thinspace \ boundary cases.

\section{An abstract finite time extinction criterion\label%
{S.general.criterion}}

The purpose of this section is to establish a general sufficient criterion
for extinction in finite time for a SBM $\,X$\thinspace \ in $\,{\sf R}^d$%
\thinspace \ with branching rate functional $\,K$.\thinspace $\;$The central
idea is to divide a finite time interval into an infinite number of stages
in such a way that all of the mass will be dead at the end of all these
stages. For this purpose, at each stage we distinguish between good and bad
historical paths. The good paths accumulate a ``significant'' rate of
branching, so that they die by the next stage, with high probability. The
remaining bad paths may not die, but by assumption they carry a small
expected mass at the beginning of the next stage.

\subsection{Refinement: historical SBM $\widetilde{X}$\label{SS.refinement}}

To realize this concept, we have to pass to a {\em ``historical''}\/
setting. \vspace{1pt}That is, the measures $\,X_t({\rm d}b)$\thinspace \ on $%
\,{\sf R}^d$\thinspace \ are thought of to be projections of measures $\,%
\widetilde{X}_t({\rm d}w)$\thinspace \ where $\,w$\thinspace \ is a
continuous function on the interval $\,[0,t].$\thinspace \ Heuristically, a
particle in $\,X_t$\thinspace \ with position $\,b$\thinspace \ is
additionally equipped with a path $\,w:[0,t]\rightarrow {\sf R}^d$\thinspace
\ with terminal point $\,w_t=b,$\thinspace \ which gives the spatial past
history of the particle and its ancestors. (For a more detailed exposition,
we refer e.g.\ to \cite{FleischmannMueller97.infcatmass}.)\smallskip 

Equip $\,{\bf C}:={\cal C}[{\sf R}_{+\,},{\sf R}^d]$\thinspace \ with the
topology of uniform convergence on all compact subsets of $\,{\sf R}_{+\,}.$%
\thinspace \ To each $\,w\in {\bf C}$\thinspace \ and $\,t\ge 0,$\thinspace
\ we associate the {\em stopped path}\/ $\,w^t\in {\bf C}$ \thinspace by
setting $\,w_s^t:=w_{t\wedge s\,},$ $\,s\ge 0.$\thinspace \ \vspace{1pt}%
Write $\,{\bf C}^t$\thinspace \ for the closed subspace of $\,{\bf C}$%
\thinspace \ of all these paths stopped at time \thinspace $t.$\thinspace \
Note that $\,{\bf C}^t$\thinspace \ could be regarded as $\,{\cal C}\big[%
[0,t],{\sf R}^d\big]$\thinspace \ (as we did in the previous paragraph), and 
$\,{\bf C}^0$\thinspace \ as $\,{\sf R}^d.$

To every $\,w\in {\bf C}$\thinspace \ we associate the corresponding {\em %
stopped path trajectory \thinspace }$\widetilde{w}$\thinspace \ by setting $%
\,\widetilde{w}_t:=w^t,$ $\,t\ge 0.$\thinspace \ Writing $\,\Vert \cdot
\Vert _\infty $\thinspace \ for the supremum norm, for $\,0\le s\le t$%
\thinspace $\;$we get 
\[
\Vert \widetilde{w}_t-\widetilde{w}_s\Vert _\infty =\Vert w^t-w^s\Vert
_\infty =\sup\limits_{s\le r\le t}\,|w_r-w_s|\longrightarrow 0\quad {\rm as}%
\quad t-s\downarrow 0.
\]
Hence, $\,\widetilde{w}\,\,\;$belongs to $\,{\cal C}[{\sf R}_{+\,},{\bf C}].$

The image of the Brownian motion $\,W$\thinspace \ under the map $\,w\mapsto 
\widetilde{w}\,\,\;$is called the {\em Brownian path process} 
\[
\widetilde{W}=\left[ \widetilde{W},\widetilde{\Pi }_{s,w\,},s\ge 0,\;w\in 
{\bf C}^s\right] .
\]
That is, at time $\,s$\thinspace \ we start with a \vspace{1pt}path $\,w=%
\widetilde{W}_s$\thinspace \ stopped at time $\,s,\,\,\;$and let a path
trajectory $\,\big\{\widetilde{W}_t:t\ge s\big\}$\thinspace \ evolve with
law $\,\widetilde{\Pi }_{s,w}\,\,\;$determined by the path $\,\{W_t:\,t\ge
s\}$\thinspace \ starting at time $\,s$\thinspace \ from $\,w_{s\,}.$

Note that if $\,K$\thinspace \ belongs to $\,{\bf K}^\nu $\thinspace \ for
some $\,\nu >0,$\thinspace \ then $\,K$\thinspace \ is also a continuous
additive functional with respect to the Brownian path \vspace{1pt}process $\,%
\widetilde{W}.$

Set 
\begin{equation}
{\sf R}_{+}\widehat{\times }{\bf C}^{\bullet }\,:=\,\left\{
(s,w):_{\!_{\!_{\,_{}}}}\,s\in {\sf R}_{+\,},\,\;w\in {\bf C}^s\right\} 
\label{cross.hut}
\end{equation}
and write $\,{\cal M}\!\left[ {\sf R}_{+}\widehat{\times }{\bf C}^{\bullet
}\right] $\thinspace \ for the set of all measures $\,\eta $\thinspace \ on $%
\,\!{\sf R}_{+}\widehat{\times }{\bf C}^{\bullet }$\thinspace \ which are
finite if restricted to a finite time interval. Moreover, let 
\begin{equation}
\widetilde{\Pi }_\eta \,:=\,\int \!\!\eta ({\rm d}s,{\rm d}w)\,\widetilde{%
\Pi }_{s,w\,},\qquad s\ge 0,\quad \eta \in {\cal M}\!\left[ {\sf R}_{+}%
\widehat{\times }{\bf C}^{\bullet }\right] ,  \label{Pi.tilde.s.mu}
\end{equation}
and 
\begin{equation}
\widetilde{\Pi }_{s,\mu }\,:=\,\widetilde{\Pi }_{\delta _s\times \mu
\,},\qquad s\ge 0,\quad \mu \in {\cal M}_{{\rm f}}[{\bf C}^s].
\end{equation}
$W$\thinspace \ can be reconstructed from $\,\widetilde{W}$\thinspace \ by 
{\em projection:}\/ $\,W_t:=(\widetilde{W}_t)_{t\,}.$\thinspace \ This will
often be used in the sequel.

Now we give the following historical version of Lemma \ref{L.cont.SBM.stable}%
, which follows from a modification of Propositions 1 (p.225), 12 (p.230),
and Lemma 4 (p.232) in \cite{DawsonFleischmann97.persist}.

\begin{proposition}[historical SBM with branching rate functional $K$]
\label{P.cont.historical.SBM}\ \hfill \newline
Let\/ $\,d\ge 1,$\thinspace \ and fix $\,K\in {\bf K}^\nu $\thinspace \ for
some $\,\nu >0.$\vspace{-5pt}

\begin{description}
\item[(a) (existence) ]  There exists a (time-inhomogeneous) Markov process$%
\,$%
\[
\widetilde{X}=\Big[\widetilde{X},\,\widetilde{P}_{s,\mu \,},\,s\ge 0,\,\mu
\in {\cal M}_{{\rm f}}[{\bf C}^s]\Big]
\]
with states\/ $\,\widetilde{X}_t\in {\cal M}_{{\rm f}}[{\bf C}^t]$, $\,t\ge
s,$\thinspace \ and with Laplace functional 
\begin{equation}
\widetilde{P}_{s,\mu }\exp 
\big\langle %
\widetilde{X}_{t\,},-\varphi 
\big\rangle
                                                                                                                                         %
\,=\,\,\exp \left\langle \mu \,,_{\!_{\!_{\,}}}-v(s,\cdot
\,|\,t)\right\rangle ,  \label{log.Laplace.X.hist}
\end{equation}
$0\le s\le t,\,\;\mu \in {\cal M}_{{\rm f}}[{\bf C}^s],\,\;\varphi \in b%
{\cal B}_{+}[{\bf C}],\,\;$where for\/ $\,t,\varphi $\thinspace \ fixed,\/ $%
\,v=v(\cdot \,,\cdot \,|\,t)\ge 0$\thinspace \ is uniquely determined by the
log-Laplace equation 
\begin{equation}
v(s,\omega _s)\,=\,\widetilde{\Pi }_{s,\omega _s}\Big[\varphi (\widetilde{W}%
_t)-\int_s^t\!K({\rm d}r)\,v^2(r,\widetilde{W}_r)\Big] \,,  \label{equ.hist}
\end{equation}
$0\le s\le t,\,\;\omega _s\in {\bf C}^s.$\vspace{-3pt}

\item[(b) (modification) ]  To each\/ $\,\eta \in {\cal M}\!\left[ {\sf R}%
_{+}\widehat{\times }{\bf C}^{\bullet }\right] $\thinspace \ there is a
Markov process\/ $\,\left[ \widetilde{X},\widetilde{P}_\eta \right] $%
\thinspace \ with states\/ $\,\widetilde{X}_t\in {\cal M}_{{\rm f}}[{\bf C}%
^t]$\thinspace \ and such that 
\begin{equation}
\widetilde{P}_\eta \exp 
\big\langle %
\widetilde{X}_{t\,},-\varphi 
\big\rangle
                                                                                                                                         %
\,=\,\,\exp \left\langle \eta \,,_{\!_{\!_{\,}}}-v(\cdot ,\cdot
\,|\,t)\right\rangle ,\qquad t\ge 0,
\end{equation}
with\/ $\,v(s,\cdot \,|\,t)$\thinspace \ from\/ {\em (a)} if\/ $\,0\le s\le
t,$\thinspace \ and\/ $\,v(s,\cdot \,|\,t)=0$\thinspace \ otherwise.

\item[(c) (moments) ]  $\big(\widetilde{X},\widetilde{P}_{s,\mu }\big)$ 
\vspace{1pt}has finite moments of all orders. In particular, for\/ $\,\eta
\in {\cal M}\!\left[ {\sf R}_{+}\widehat{\times }{\bf C}^{\bullet }\right] $%
\thinspace \ and\/ $\,\varphi _1,\varphi _2\in b{\cal B}_{+}[{\bf C}],$%
\thinspace \ as well as\/ $\;t_1,t_2\ge 0,$%
\begin{equation}
\widetilde{P}_\eta 
\big\langle %
\widetilde{X}_{t_1\,},\varphi _1%
\big\rangle
                                                                                                                                         %
\,=\,\widetilde{\Pi }_\eta \,\varphi _1\big(\widetilde{W}_{t_1}\big)\,,
\label{tilde.expectation}
\end{equation}
\begin{equation}
\left. 
\begin{array}{l}
\widetilde{C}{\rm ov}_\eta \left[ \!%
\big\langle %
\widetilde{X}_{t_1\,},\varphi _1%
\big\rangle
                                                                                                                                         %
,%
\big\langle %
\widetilde{X}_{t_2\,},\varphi _2%
\big\rangle
                                                                                                                                         %
\!\right] \vspace{5pt} \\ 
=\,2\,\widetilde{\Pi }_\eta \dint \!\!K({\rm d}r)\Big[\,\widetilde{\Pi }_{r,%
\widetilde{W}_r}\varphi _1\big(\widetilde{W}_{t_1}\big)\Big] \Big[\,%
\widetilde{\Pi }_{r,\widetilde{W}_r}\varphi _2\big(\widetilde{W}_{t_2}\big)%
\Big] .
\end{array}
\;\right\} 
\end{equation}
\end{description}
\end{proposition}

We call this superprocess $\,\widetilde{X}$\thinspace \ the {\em historical
SBM with branching rate functional }$\,K.$

Of course, $\,X$\thinspace \ can be gained back from $\,\widetilde{X}$%
\thinspace \ by {\em projection:} 
\begin{equation}
X_t\,=\,\widetilde{X}_t\Big(\big\{w\in {\bf C}^t:\;w_t\in (\cdot )\big\}\Big)%
.  \label{projection}
\end{equation}

\subsection{Dynkin's ``stopped'' measures $\widetilde{X}_\tau $\label%
{SS.stopped.measures}}

We also have to recall Dynkin's \cite{Dynkin91.a,Dynkin91.p} concept of {\em %
``stopped''}\/ historical superprocesses. We have two reasons for this.
First, to handle also the lattice Model 2, we have to allow the times $%
\,T_1<T_2<\cdots $\thinspace \ mentioned in Subsection \ref{SS.Model1} to be
Brownian stopping times. The second reason is that we intend to scale the
SBM along individual particles' trajectories according to their accumulated
rate of branching.

Roughly speaking, if $\,\tau $\thinspace \ is a {\em Brownian}\/ stopping
time, Dynkin's stopped measure $\,\widetilde{X}_\tau $\thinspace \ describes
the population one gets, if each (individual) reactant path is stopped in
the moment $\,\tau .$\thinspace \ For a detailed development, we refer to 
\cite[Subsection 1.5]{Dynkin91.a} and \cite[Subsection 1.10]{Dynkin91.p}.
For convenience, here we collect only the following facts.

Let $\,\tau _{t\,},$ $\,t\ge 0,$\thinspace \ be {\em stopping times}\/ with
respect to the (natural) filtration of Brownian motion $\,W,$\thinspace \
satisfying $\,\tau _s\le \tau _t$\thinspace \ if $\;s\le t.$\thinspace \
Then there is a family 
\[
\left\{ \widetilde{X}_{\tau _t}:\,\,t\ge 0\right\} ,
\]
of random measures in $\,{\cal M}\!\left[ {\sf R}_{+}\widehat{\times }{\bf C}%
^{\bullet }\right] ,\,\,\;$the so-called \label{warning}{\em ``stopped''
historical SBM}\/ related to the family of Brownian stopping times
\thinspace $\left\{ \tau _t:\,t\ge 0\right\} .$\thinspace \ This family
satisfies the so-called {\em special Markov property,}\/ which roughly says
the following. For $\,\vspace{2pt}s\ge 0,$\thinspace \ let $\,{\cal G}_{\tau
_s}$\thinspace \ denote the pre-$\tau _s$\thinspace \ $\sigma $--field
(concerning the historical superprocess $\,\widetilde{X}).$\thinspace \
Given $\,{\cal G}_{\tau _{s\,}},$\thinspace \ hence in particular $\,%
\widetilde{X}_{\tau _s}=:\vartheta ,$\thinspace \ the stopped process $\,%
\big\{\widetilde{X}_{\tau _t}:\,\,t\ge s\big\}
$\thinspace \ starts anew \vspace{1pt}(\cite[Theorem 1.6]{Dynkin91.a} and 
\cite[Theorem 1.5]{Dynkin91.p}), namely based on the law $\,\widetilde{P}%
_{\vartheta \,}.$\vspace{1pt}

Similarly, the notation of a sequence $\,\big\{\widetilde{X}_{\tau
_n}:\,\,n\ge 1\big\}
$\thinspace \ of stopped measures related to Brownian stopping times $\,\tau
_1\le \tau _2\le \cdot \cdot \cdot \;$\ can be introduced.

In formal analogy with Proposition \ref{P.cont.historical.SBM}\thinspace
(c), the following first two {\em moment formulas}\/ hold (\cite[(1.50a)]
{Dynkin91.a}). For $\,\,t\ge 0$\thinspace \ and $\,\eta \in {\cal M}\!\left[ 
{\sf R}_{+}\widehat{\times }{\bf C}^{\bullet }\right] $\thinspace \ as well
as $\,\varphi $\thinspace \ in $\,b{\cal B}_{+}[{\bf C}],$ 
\begin{eqnarray}
\widetilde{P}_\eta 
\big\langle %
\widetilde{X}_{\tau _t\,},\varphi 
\big\rangle
                                                                                                                                         %
&=&\widetilde{\Pi }_\eta \,\varphi \big(\,\widetilde{W}_{\tau _t}\big)\,,
\label{stopped.expectation} \\[2pt]
\widetilde{V}\!{\rm ar}_\eta 
\big\langle %
\widetilde{X}_{\tau _t\,},\varphi 
\big\rangle
                                                                                                                                         %
&=&2\,\widetilde{\Pi }_\eta \dint \!\!K({\rm d}r)\Big[\,\widetilde{\Pi }_{r,%
\widetilde{W}_r}\varphi \big(\,\widetilde{W}_{\tau _t}\big)\Big] ^2.
\label{stopped.variance}
\end{eqnarray}

\subsection{The method of good and bad historical paths\label%
{SS.method.paths}}

Fix $\,K\in {\bf K}^\nu ,$\thinspace \ for some $\,\nu >0,$\thinspace \ and
a finite measure $\,\mu $\thinspace \ on $\,{\sf R}^d.$\thinspace \ Consider
the historical SBM $\,\widetilde{X}$\thinspace \ of Proposition \ref
{P.cont.historical.SBM} starting from $\,\widetilde{X}_0=\mu .\,\;$First we
introduce some Brownian stopping times and small constants.

\begin{hypothesis}[stage quantities]
\label{H.stage.quantities}{\em Let} $\,$\thinspace $0<\varepsilon <1$%
\thinspace \ {\em and }$\,N=N(\varepsilon )\geq 0.$\vspace{-6pt}\nolinebreak 

\begin{description}
\item[(a) (stage duration) ]  {\em Consider} Brownian stopping times\/ $\,$%
\begin{equation}
0\le T_N^\varepsilon <T_{N+1}^\varepsilon <\cdots <T_\infty ^\varepsilon
<\infty ,  \label{stopping.epsilon}
\end{equation}
{\em where the bound }$\,T_\infty ^\varepsilon $\thinspace \ {\em is} non-%
{\em random.}$\vspace{-4pt}$

\item[(b) (constants) ]  {\em For }$\,n\ge N=N(\varepsilon ),$\thinspace \ 
{\em let }$\,M_{n\,}^\varepsilon ,\,\delta _{n\,}^\varepsilon ,\,\lambda
_{n\,}^\varepsilon ,\,\xi _{n\,}^\varepsilon \;{\em and,}${\em \ in addition,%
}\/ $\,\delta _{N-1}^\varepsilon $\thinspace \ {\em be (strictly) positive
constants with the following properties:}\vspace{-4pt}

\begin{description}
\item[(b1) ]  \ \hfill $M_n^\varepsilon \downarrow 0\quad {\rm as}\quad
n\uparrow \infty .\;$\hfill \ 

\item[(b2) ]  \ \hfill $\lim\limits_{\varepsilon \downarrow 0}\,\bigg(\delta
_{N-1}^\varepsilon +\dsum\limits_{n\ge N}\left( \delta _n^\varepsilon
+\lambda _n^\varepsilon \right) \bigg)\,=\,0.$\hfill $\Diamond $\smallskip 
\end{description}
\end{description}
\end{hypothesis}

Introduce the set $\,E_n^\varepsilon $\thinspace \ of so-called {\em good
historical paths }($\!${\em during }$\,[T_{n\,}^\varepsilon
,T_{n+1}^\varepsilon ]),$ 
\begin{equation}
E_n^\varepsilon \,:=\,\bigg\{w\in {\bf C:\;}\int_{T_n^\varepsilon
}^{T_{n+1}^\varepsilon }K({\rm d}r)\ge \xi _n^\varepsilon \bigg\},
\label{E.n}
\end{equation}
that is, paths with at least the amount $\,\xi _n^\varepsilon $\thinspace \
of accumulated branching over the time interval $\,[T_n^\varepsilon
,T_{n+1}^\varepsilon ).$\thinspace \ We call $\,(E_n^\varepsilon )^{{\rm c}}=%
{\bf C}\backslash E_n^\varepsilon $\thinspace \ the set of {\em bad}\/
paths. On the good and bad paths we impose the following hypothesis.

\begin{hypothesis}[good and bad paths]
\label{H.good.and.bad}{\em Fix }$\,\varepsilon \in (0,1).$\thinspace \ {\em %
First of all,} 
\begin{equation}
\widetilde{P}_{0,\mu }\left( \big\| \widetilde{X}_{T_N^\varepsilon }\big\| %
\,>\,M_N^\varepsilon \right) \;\le \;\delta _{N-1\,}^\varepsilon ,
\label{starting.cond}
\end{equation}
{\em and for all }$\,n\ge N=N(\varepsilon ),$%
\begin{equation}
\widetilde{P}_{0,\mu }\left\{ \widetilde{X}_{T_{n+1}^\varepsilon
}(E_n^\varepsilon )>0\;\Big| \;\big\| \widetilde{X}_{T_n^\varepsilon }\big\| %
\,\le \,M_n^\varepsilon \right\} \;\le \;\delta _{n\,}^\varepsilon ,
\label{Feller.survival}
\end{equation}
\begin{equation}
\widetilde{P}_{0,\mu }\left\{ \widetilde{X}_{T_{n+1}^\varepsilon }\left(
(E_n^\varepsilon )_{\!_{\!_{\,}}}^{{\rm c}}\right) \;\Big| \;\big\| 
\widetilde{X}_{T_n^\varepsilon }\big\| \,\le \,M_n^\varepsilon \right\}
\;\le \;\lambda _n^\varepsilon \,M_{n+1}^\varepsilon \,.\vspace{-10pt}
\label{small.expectation}
\end{equation}
\hfill 

\hfill $\Diamond $\smallskip 
\end{hypothesis}

Here is our {\em interpretation}\/ of Hypothesis \ref{H.good.and.bad}.
Recall that by Hypothesis \ref{H.stage.quantities}\thinspace (b) the \vspace{%
1pt}numbers $\,M_{n\,}^\varepsilon ,\,\delta _{N-1\,}^\varepsilon ,\,\delta
_{n\,}^\varepsilon ,$\thinspace \ and $\,\lambda _n^\varepsilon $\thinspace
\ are small. So at the beginning of the $\,N^{{\rm th\,}}$ stage the total
mass $\,\big\| 
\widetilde{X}_{T_N^\varepsilon }\big\| \,\;$is already small with a high $\,%
\widetilde{P}_{0,\mu }$--probability. Then starting with a small mass at the
beginning of the $\,n^{{\rm th\,}}$ stage, our condition (\ref
{Feller.survival}) says that {\em good paths survive}\/ only with a small
(conditional) probability in the present stage, whereas (\ref
{small.expectation}) means that the (conditional) {\em expected mass of bad
paths}\/ is small.\vspace{2pt}

Our abstract criterion now reads as follows. Recall that $\,d=1,$\thinspace
\ $\mu \in {\cal M}_{{\rm f\,}},$\thinspace \ and that $\,K\in {\bf K}^\nu $%
\thinspace \ for some $\,\nu >0.$

\begin{theorem}[abstract criterion for finite time extinction]
\label{T.abstract.criterion}Impose Hy{\Large \-}potheses {\em \ref
{H.stage.quantities}} and\/ {\em \ref{H.good.and.bad}}. Then with\/ $\,%
\widetilde{P}_{0,\mu }$--probability one,\/ $\,\widetilde{X}_t\,=0\,\;\,$for
some\/ $\,t.$
\end{theorem}

We mention that under additional conditions, estimate (\ref{Feller.survival}%
) can be obtained by a comparison with Feller's branching diffusion, see
Subsection \ref{SS.comparison} below, whereas the expectation formula for
stopped historical SBM is available to reduce assertion (\ref
{small.expectation}) to a statement on the probability of a path to be bad,
i.e.\ to have a small accumulated rate of branching. In fact, by the special
Markov property and the expectation formula (\ref{stopped.expectation})
applied to the indicator function $\,\varphi ={\bf 1}_{(E_n^\varepsilon )^{%
{\rm c}}}$\thinspace \ and the starting measure $\,\eta =\delta
_{T_n^\varepsilon }\times \widetilde{X}_{T_n^\varepsilon \,},$\thinspace \
we have 
\[
\widetilde{P}_{0,\mu \,}\left\{ \widetilde{X}_{T_{n+1}^\varepsilon }\left(
(E_n^\varepsilon )_{\!_{\!}}^{{\rm c}}\right) \;\Big| \;{\cal G}%
_{T_n^\varepsilon }\right\} \,\,=\,\int \!\!\widetilde{X}_{T_n^\varepsilon }(%
{\rm d}w)\;\widetilde{\Pi }_{T_n^\varepsilon \,,w}\left( \widetilde{W}%
_{T_{n+1}^\varepsilon }\in _{\!_{\!_{\,}}}(E_n^\varepsilon )^{{\rm c}%
}\right) .
\]
Since $\,E_n^\varepsilon $\thinspace \ only depends on $\,\left\{ w_s:\,s\ge
T_n^\varepsilon \right\} ,$\thinspace \ we can write 
\begin{equation}
\widetilde{P}_{0,\mu \,}\left\{ \widetilde{X}_{T_{n+1}^\varepsilon }\left(
(E_n^\varepsilon )_{\!_{\!}}^{{\rm c}}\right) \;\Big| \;{\cal G}%
_{T_n^\varepsilon }\right\} \,\,=\,\int \!\!X_{T_n^\varepsilon }({\rm d}%
a)\;\Pi _{T_n^\varepsilon ,a}\left( W\in _{\!_{\!_{\,}}}(E_n^\varepsilon )^{%
{\rm c}}\right)   \label{conditional.expectation}
\end{equation}
(recall Remark \ref{R.inhomogeneous}). Then (\ref{conditional.expectation})
implies the following result.

\begin{lemma}[sufficient condition]
\label{L.new}If the estimate 
\begin{equation}
\Pi _{T_n^\varepsilon ,a}\left( W\in _{\!_{\!_{\,}}}(E_n^\varepsilon )^{{\rm %
c}}\right) \;\le \;\lambda _n\,\frac{M_{n+1}^\varepsilon }{M_n^\varepsilon }%
,\qquad a\in {\rm supp}X_{T_n^\varepsilon \,},  \label{Carl's.addition}
\end{equation}
holds, then the conditional expectation estimate\/ {\em (\ref
{small.expectation})} is true.
\end{lemma}

\subsection{Proof of the abstract criterion\label{SS.proof.abstract}}

Here we want to prove Theorem \ref{T.abstract.criterion}. For $%
\,0<\varepsilon <1$\thinspace \ and $\,N=N(\varepsilon )\ge 0,$\thinspace \
set 
\begin{equation}
A_n^\varepsilon \,:=\,\left\{ \big\| \widetilde{X}_{T_n^\varepsilon }\big\| %
\le M_n^\varepsilon \right\} ,\quad n\ge N,\quad {\rm and}\quad
A^\varepsilon \,:=\,\bigcap_{n\ge N}A_{n\,}^\varepsilon ,
\end{equation}
as well as 
\begin{equation}
\overline{T}_\infty ^{\,\varepsilon }\,:=\,\lim_{n\uparrow \infty
}T_{n\,}^\varepsilon .
\end{equation}
Note that this limiting Brownian stopping time satisfies $\,\overline{T}%
_\infty ^{\,\varepsilon }\le T_\infty ^\varepsilon <\infty $.\medskip 

\noindent 1$^{\circ }$ ({\em extinction on}\/ $\,A^\varepsilon $)\quad First
of all we show that for all $\,\varepsilon \in (0,1),$%
\begin{equation}
\big\| \widetilde{X}_{\overline{T}_\infty ^{\,\varepsilon }}\big\| %
\,=\,0\quad {\rm on\;\,}A^\varepsilon ,\quad \widetilde{P}_{0,\mu }{\rm -a.s.%
}  \label{onA}
\end{equation}
Indeed, fix $\,0<\varepsilon <1$\thinspace \ and a $\,\zeta >0.$\thinspace \
Then by Markov's inequality, for each $\,n\ge N,$%
\begin{equation}
\widetilde{P}_{0,\mu }\left( \left\{ \big\| \widetilde{X}_{\overline{T}%
_\infty ^{\,\varepsilon }}\big\| \,>\zeta \right\} \cap A^{\varepsilon
^{\!^{\!^{\,^{\!}}}}}\right) \;\le \;\zeta ^{-1}\,\widetilde{P}_{0,\mu \,}%
{\bf 1}_{A_{n\,}^\varepsilon }\big\| \widetilde{X}_{\overline{T}_\infty
^{\,\varepsilon }}\big\| .  \label{Markov}
\end{equation}
But by the special Markov property, the expectation formula (\ref
{stopped.expectation}), and the definition of $\,A_{n\,}^\varepsilon ,$%
\[
\widetilde{P}_{0,\mu }\left\{ {\bf 1}_{A_{n\,}^\varepsilon }\big\| 
\widetilde{X}_{\overline{T}_\infty ^{\,\varepsilon }}\big\| \;\Big| \;{\cal G%
}_{T_n^\varepsilon }\right\} \;=\;{\bf 1}_{A_n^\varepsilon \,}\widetilde{P}_{%
\widetilde{X}_{T_n^\varepsilon }}\big\| 
\widetilde{X}_{\overline{T}_\infty ^{\,\varepsilon }}\big\| \;=\;{\bf 1}%
_{A_{n\,}^\varepsilon }\big\| 
\widetilde{X}_{T_n^\varepsilon }\big\| \;\le \;M_n^\varepsilon \,.
\]
Hence, estimate (\ref{Markov}) can be continued with 
\[
\le \;\zeta ^{-1}M_n^\varepsilon \;\stackunder{n\uparrow \infty }{%
\longrightarrow }\;0,
\]
by Hypothesis \ref{H.stage.quantities}\thinspace (b1). Thus, 
\[
\widetilde{P}_{0,\mu }\left( \left\{ \big\| \widetilde{X}_{\overline{T}%
_\infty ^{\,\varepsilon }}\big\| \,>\zeta \right\} \cap A^{\varepsilon
^{\!^{\!^{\,^{\!}}}}}\right) \,=\,0\quad \forall \,\zeta >0,
\]
and (\ref{onA}) follows.\medskip 

\noindent 2$^{\circ }$ ($A_{n+1}^\varepsilon $\thinspace \ \vspace{1pt}{\em %
fails with small conditional probability})\quad From Markov's inequality and
(\ref{small.expectation}) we get, for $\,0<\varepsilon <1$\thinspace \ and $%
\,n\ge N,$%
\[
\widetilde{P}_{0,\mu }\left\{ \widetilde{X}_{T_{n+1}^\varepsilon }\left(
(E_n^\varepsilon )_{\!_{\!_{\,}}}^{{\rm c}}\right) \,>\,M_{n+1}^\varepsilon
\;\Big| \;A_n^\varepsilon \right\} \;\le \;{\lambda }_{n\,}^\varepsilon .
\]
Together with (\ref{Feller.survival}), we conclude for 
\[
\widetilde{P}_{0,\mu }\left\{ \left( A_{n+1}^\varepsilon \right) ^{{\rm c}%
^{\!^{}}}\;\big| \;A_n^\varepsilon \right\} \;\le \;\delta _n^\varepsilon +{%
\lambda }_n^\varepsilon \,.\vspace{4pt}
\]

\noindent 3$^{\circ }$ ($A^\varepsilon $\thinspace \ {\em fails with small
probability})\quad Next we show that 
\begin{equation}
\lim_{\varepsilon \downarrow 0}\,\widetilde{P}_{0,\mu }\left( \left(
A^\varepsilon \right) ^{{\rm c}^{\!^{}}}\right) \,=\,0.  \label{largeA}
\end{equation}
We decompose the complement $\,\left( A^\varepsilon \right) ^{{\rm c}}$%
\thinspace \ of $\,A^\varepsilon $\thinspace \ according to the smallest
natural number $\,n\ge N$\thinspace \ such that $\,A_n^\varepsilon $%
\thinspace \ fails: 
\begin{eqnarray*}
\widetilde{P}_{0,\mu }\left( \left( A^\varepsilon \right) ^{{\rm c}}\right) 
&=&\widetilde{P}_{0,\mu }\left( \left( A_N^\varepsilon \right) ^{{\rm c}%
^{\!^{}}}\;\right) +\sum_{n\ge N}\,\widetilde{P}_{0,\mu }\left(
A_N^\varepsilon \cap \cdot \cdot \cdot \cap A_n^\varepsilon \cap \left(
A_{n+1}^\varepsilon \right) ^{{\rm c}^{\!^{}}}\right)  \\
&\le &\widetilde{P}_{0,\mu }\left( \left( A_N^\varepsilon \right) ^{{\rm c}%
^{\!^{}}}\;\right) +\sum_{n\ge N}\,\widetilde{P}_{0,\mu }\left\{ \left(
A_{n+1}^\varepsilon \right) ^{{\rm c}^{\!^{}}}\;\big| \;A_n^\varepsilon
\right\} .
\end{eqnarray*}
Then (\ref{largeA}) follows from (\ref{starting.cond}), step 2$^{\circ }$
and Hypothesis \ref{H.stage.quantities}\thinspace (b2).\medskip 

\noindent 4$^{\circ }$ ({\em conclusion})\quad Let $\,\varepsilon =k^{-1},$%
\thinspace \ $k>1.$\thinspace \ From (\ref{largeA}) and the monotonicity
property of \vspace{1pt}measures, we learn that the event that $\,A^{1/k}$%
\thinspace \ fails for all $\,k,$\thinspace \ has $\widetilde{P}_{0,\mu }$%
--probability $\,0.$\thinspace \ In other words, 
\[
\widetilde{P}_{0,\mu }\bigg(\bigcup_{k>1}A^{1/k}\bigg)\,=\,1.
\]
Then step 1$^{\circ }$ implies that 
\[
\exists \;k>1\;\,{\rm such\;that\;\,}\widetilde{X}_{\overline{T}_\infty
^{\,1/k}}=0,\quad \widetilde{P}_{0,\mu }{\rm -a.s.}
\]
Again applying the special Markov property, we obtain 
\[
\exists \;k\;\,{\rm such\;that\;\,}\widetilde{X}_{T_\infty ^{\,1/k}}=0,\quad 
\widetilde{P}_{0,\mu }{\rm -a.s.}
\]
Since $\,T_\infty ^{\,1/k}$\thinspace \ is non-random, the proof of our
abstract Theorem \ref{T.abstract.criterion} is finished.\hfill 
\endproof%

\subsection{Scaled comparison with Feller's branching diffusion\label%
{SS.comparison}}

Consider a pair of Brownian stopping times $\,0\le T_0<T_1,$\thinspace \ a
constant $\,\xi >0,$\thinspace \ and define $\,$%
\[
E\,:=\,\bigg\{w\in {\bf C:\;}\int_{T_0}^{T_1}K({\rm d}r)\ge \xi \bigg\}.
\]
We want to estimate a conditional probability as in (\ref{Feller.survival})
in Hypothesis \ref{H.good.and.bad} under a mild additional assumption on the
branching rate functional $\,K.$\thinspace \ For this purpose, we will
compare with the survival probability in Feller's branching diffusion.

Recall that $\,d\ge 1,$\thinspace \ $\mu \in {\cal M}_{{\rm f\,}},$%
\thinspace \ and that $\,{\cal G}_{T_0\,}\;$denotes the pre-$T_0$\thinspace
\ $\sigma $--field.

\begin{proposition}[comparison with Feller's branching diffusion]
\label{P.comparison.Feller}\hspace{-2.7pt}Assume that the branching rate
functional\/ $\,K\in {\bf K}^\nu $\thinspace \ $(\nu >0)$\thinspace \ is
homogeneous and satisfies 
\begin{equation}
\int_0^t\!K({\rm d}r)\stackunder{t\uparrow \infty }{\longrightarrow }\infty
,\quad \Pi _{0,a}{\rm -}a.s.,\quad a\in {\sf R}^d{\sf .}
\label{infinite.BCLT}
\end{equation}
Then $\,\widetilde{P}_{0,\mu }$--almost surely, 
\begin{equation}
\widetilde{P}_{0,\mu }\left\{ \widetilde{X}_{T_1}(E)>0\;\Big| \;{\cal G}%
_{T_0}\right\} \,\le \,\,\frac 1\xi \,\Vert \widetilde{X}_{T_0}\Vert \,.%
\vspace{2pt}  \label{3.5}
\end{equation}
\end{proposition}

We will prove Proposition \ref{P.comparison.Feller} in the next subsection,
using an idea from \cite{FleischmannMueller97.infcatmass}, which was in turn
inspired by a modulus of continuity technique of \cite{DawsonPerkins91}. In
fact, since the paths in $\,E$\thinspace \ have a ``significant''
accumulated rate of branching over the time interval $\,[T_0,T_1)$\thinspace
\ [recall (\ref{E.n})], we can {\em compare}\/ (in law) $\,\widetilde{X}%
_{T_1}(E)$\thinspace \ with the mass in Feller's branching diffusion after
an appropriate individual time change. (Recall that the total mass of the
classical super-Brownian motion is equal in distribution to Feller's
branching diffusion.) But for Feller's branching diffusion, there is a
well-known estimate for the probability that the process survives in the
given time.

\begin{remark}
{\em We usually apply Proposition \ref{P.comparison.Feller} for }$%
\,T_0=T_n^\varepsilon ${\em \thinspace }$,${\em \ }$\,T_1=T_{n+1\,}^%
\varepsilon ${\em \thinspace ,\ }$\xi =\xi _{n\,}^\varepsilon ,${\em %
\thinspace \ and }$\,E=E_{n\,}^\varepsilon ,${\em \thinspace \ with fixed }$%
\,\varepsilon $\thinspace \ {\em and }$\,n.$\thinspace \ \hfill $\Diamond $
\end{remark}

\subsection{Proof of the comparison argument\label{SS.proof.comparison}}

Consider the process $\,\big\{\widetilde{X}_t:\,t\ge T_0\big\},$\thinspace \
given $\,{\cal G}_{T_{0\,}}.$\thinspace \ In particular, the starting
measure $\,\widetilde{X}_{T_0}=:\vartheta $\thinspace \ is given. Note that
this (conditional) process has the law $\,\widetilde{P}_{\vartheta \,},$%
\thinspace \ by the special Markov property. In order to prove Proposition 
\ref{P.comparison.Feller}, we first intend to define a {\em new time scale}%
\/ denoted by $\,r,\,\,\;$dictated by the additive functional $\,K.$%
\thinspace \ Given for the moment a path $\,w\in {\bf C},$\thinspace \ set%
\vspace{-2pt} 
\begin{equation}
R(t):=\int_{T_0}^{T_0+t}K({\rm d}s),\qquad t\ge 0,  \label{R}
\end{equation}
(recall that $\,K$\thinspace \ is a continuous additive functional of
Brownian motion $\,W).$ Note that $\,R(t)\uparrow \infty \,$\thinspace \ as $%
\,\,t\uparrow \infty ,$\thinspace \ $\widetilde{\Pi }_\vartheta $--almost
everywhere [by assumption (\ref{infinite.BCLT})], and that $\,R(t)$%
\thinspace \ depends continuously on $\,t\,\,\;[$by the continuity of $\,K].$%
\thinspace \ Define finite (Brownian) stopping times $\,\tau (r)$ $\;$%
(converging to infinity as $\,r\uparrow \infty )\,\,\;$by 
\begin{equation}
\tau (r)\,:=\,\inf \left\{ t>0:_{\!_{\!_{\,_{}}}}\,R(t)\ge r\right\} ,\qquad
r\ge 0.  \label{sigma}
\end{equation}
Consider the stopped historical SBM $\,r\mapsto \widetilde{X}_{T_0+\tau
(r)\,}.$\thinspace \ Put 
\begin{equation}
Z_r\,:=\,\big\| \widetilde{X}_{T_0+\tau (r)}\big\| ,\qquad r\ge 0,
\label{Zt}
\end{equation}
for its total mass process. Assume for the moment that $\widetilde{P}_{0,\mu
}$--a.s.\ under the probability laws $\,\widetilde{P}_\vartheta $\thinspace
\ the following two statements hold:\vspace{-4pt}

\begin{description}
\item[(i) ]  \ If $\,Z_\xi =0$\thinspace \ then $\,\widetilde{X}_{T_1}(E)=0$%
\thinspace \ (extinction of good paths).\vspace{-4pt}

\item[(ii) ]  The process $\,r\mapsto Z_r$\thinspace \ satisfies equation (%
\ref{Feller's.diffusion}) at p.\pageref{Feller's.diffusion} with $%
\,Z_0=\left\| \vartheta \right\| .$\vspace{-4pt}
\end{description}

\noindent Then, from the well-know survival probability formula for
solutions $\,Z$\thinspace \ of equation (\ref{Feller's.diffusion}), that is
of Feller's branching diffusion, we have 
\begin{equation}
\widetilde{P}_\vartheta \Big(\widetilde{X}_{T_1}(E)>0\Big)\,\le \,\widetilde{%
P}_\vartheta \left( Z_\xi >0\right) \,=\,\,1-\exp \left[ -\xi
^{-1^{\!^{}}}\,\left\| \vartheta \right\| \right] \;\le \;\xi ^{-1}\,\left\|
\vartheta \right\| ,  \label{Feller's.survival}
\end{equation}
which would imply (\ref{3.5}).

It remains to prove the \vspace{1pt}statements {\bf (i)} and {\bf (ii)}. To
show {\bf (i)}, assume that 
\[
\big\| \widetilde{X}_{T_0+\tau (\xi )}\big\| \,=\,Z_\xi \,=\,0.
\]
That is, all paths which accumulated the rate of branching $\,\xi $%
\thinspace \ have died [by time $\,T_0+\tau (\xi )].$\thinspace \ Therefore
we will not find paths with accumulated rate of branching greater or equal
to $\,\xi ,$\thinspace \ \vspace{1pt}in particular at time $\,T_1\ge
T_0+\tau (\xi ),$\thinspace \ which holds under $\,E.$\thinspace \
Consequently, $\,\widetilde{X}_{T_1}(E)=0,\,\;$and {\bf (i)} is verified.\ 
\vspace{1pt}

We are left with proving {\bf (ii)}. The initial condition is trivially
fulfilled. To simplify notation, we write $\,{\cal G}_r$\thinspace \ for the
pre-[$T_0+\tau (r)]$\thinspace \ $\sigma $--field. It is sufficient to show
that $\,\widetilde{P}_{0,\mu }$--almost surely, $\,Z$\thinspace \ is a $\,%
\big(\widetilde{P}_\vartheta ,\,\left( {\cal G}_r\right) _{r\ge 0}\big)${\em %
--martingale}\/ with square variation 
\[
\langle \!\langle Z\rangle \!\rangle _r\,=2\int_0^r\!{\rm d}s\;Z_s\,,\qquad
r\ge 0.
\]
This would be verified if we proved that for $\,0\le r<r^{\prime }$, 
\begin{eqnarray}
\widetilde{P}_\vartheta \left\{ Z_{r^{\prime }}\,\,\big| \,\,{\cal G}%
_r\right\}  &=&Z_{r\,},  \label{Nr1} \\[2pt]
\widetilde{P}_\vartheta \left\{ Z_{r^{\prime }}^2-2\dint_r^{r^{\prime }}\!%
{\rm d}s\;Z_s\;\bigg| \,\,{\cal G}_r\right\}  &=&Z_{r\,}^2.  \label{Nr2}
\end{eqnarray}
By the special Markov property, statement (\ref{Nr1}) follows from the
expectation formula (\ref{stopped.expectation}). But then (\ref{Nr2})
reduces to 
\begin{equation}
\widetilde{V}\!{\rm ar}_{\vartheta (r)\,}Z_{r^{\prime
}}\,\,=\,\,2\,(r-r^{\prime })\,Z_r  \label{end}
\end{equation}
with $\,\vartheta (r):=\widetilde{X}_{T_0+\tau (r)\,}.$\thinspace \ But from
the variance formula (\ref{stopped.variance}), the left hand side of (\ref
{end}) equals 
\[
2\,\,\widetilde{\Pi }_{\vartheta (r)}\int_{\tau (r)}^{\tau (r^{\prime })}\!K(%
{\rm d}s).
\]
Using the definition (\ref{sigma}) of $\,\tau (r),\,\,\;$by a change of
variables, see e.g.\ \cite[Proposition (0.4.9)]{RevuzYor91}, we can continue
with $\,\,$%
\[
=\;2\,\,\widetilde{\Pi }_{\vartheta (r)}\int_r^{r^{\prime }}\!{\rm d}%
s\;=\;2\,(r^{\prime }-r)\,\Vert \vartheta (r)\Vert ,
\]
getting the right hand side of (\ref{end}).

This completes the proof of (\ref{3.5}), that is of Proposition \ref
{P.comparison.Feller}.\hfill 
\endproof%

\section{Two applications of the abstract criterion\label{S.application}}

Here we want to apply our abstract finite time extinction criterion,
combined with the comparison with Feller's branching diffusion, to two {\em %
one-dimensional}\/ models, namely SBM with the parabolic catalyst $\,\chi _q$%
\thinspace \ of Model 3, and with a certain point catalyst $\,\underline{%
\Gamma }$\thinspace \ with atoms whose locations are dense in ${\sf R}$%
\thinspace \ (we will need $\,\underline{\Gamma }$\thinspace \ later on).

\subsection{Parabolic catalyst $\chi _q$\label{SS.application.parabolic}
(Model 3)}

In the parabolic catalyst model, the branching rate functional $\,K$%
\thinspace \ is given by the BCLT $\,L_{[W,\psi ]}\,\,\;$(recall Lemma \ref
{L.BCLT.Gamma}). Here the measure $\,\psi ({\rm d}b)\in {\cal M}_{p\,},\;$ $%
p>1,$\thinspace \ has a density function $\,\chi _q(b)=|b|^q\wedge
1,\,\;b\in {\sf R},$\thinspace \ with respect to the (normalized) Lebesgue
measure $\,{\rm d}b,$\thinspace \ as introduced in (\ref{parabolic.rate}) at
p.\pageref{parabolic.rate} (where the exponent $\,q>0$\thinspace \ is a
fixed constant).

Since the catalyst is small only near the origin $\,b=0,$\thinspace \ the 
{\em bad}\/ historical paths will be those which spend a large amount of
time near $\,0.$\thinspace \ To estimate the probability of such paths, we
need the following lemma. Let $\,L_t(b)$\thinspace \ denote the {\em local
time}\/ at $\,b\in {\sf R}$\thinspace \ of a (one-dimensional standard)
Brownian path $\,W,$\thinspace \ up to time $\,t.$

\begin{lemma}[Brownian local time large deviations]
\label{lemma3.3}There exists a const{\Large \-}ant\/ $\,c_0>0$\thinspace \
such that for all\/ $\,\theta \in (0,1],$\thinspace \ and for all\/ $\,a\in 
{\sf R},$\thinspace \ the following holds: 
\begin{equation}
\Pi _{0,a}\bigg(\int_{-\theta }^\theta {\rm d}b\;\,L_t(b)\,\ge \,\dfrac t2%
\bigg)\;\le \;\,\exp \left[ -\dfrac{c_0\,t}{\theta ^2}\right] .\vspace{2pt}
\label{varadhan}
\end{equation}
\end{lemma}

\proof%
By Brownian scaling, we have 
\[
\Pi _{0,a}\Big(\int_{-\theta }^\theta {\rm d}b\;\,L_t(b)\ge \frac t2\,\Big)%
\,=\;\Pi _{0,a/\theta }\Big(\,\int_{-1}^1{\rm d}b\;\,L_{t/\theta ^2}(b)\ge
\frac t{2\theta ^2}\,\Big).
\]
So it suffices to prove the claim for $\,\theta =1.$\thinspace \ But then we
may apply Lemma 2.2 of Donsker and Varadhan \cite{DonskerVaradhan77} with $%
\,A$\thinspace \ there the set of subprobability measures $\,\nu $\thinspace
\ on $\,{\sf R}$\thinspace \ such that $\,\nu \left(
[-1,1]_{\!_{\!_{\,}}}\right) \ge 1/2$.\hfill 
\endproof%
\bigskip 

In order to specify the quantities entering in Hypotheses \ref
{H.stage.quantities} and \ref{H.good.and.bad}, fix an $\,\alpha >0$%
\thinspace \ and a $\,\beta \in (0,2\alpha ).$\thinspace \ For $\,n\ge
N(\varepsilon )\equiv 0,$\thinspace \ and $\,0<\varepsilon <1,$\thinspace \
put 
\[
\theta _n\,:=\,{\rm e}^{-\alpha n},\quad M_n^\varepsilon \,\equiv \,M_n\,:=\,%
{\rm e}^{-(1+\beta +\alpha q)n}
\]
\[
t_n^\varepsilon \,:=\,\varepsilon ^{-1}{\rm e}^{-\beta n},\quad \xi
_n^\varepsilon \,:=\,\frac 12\,t_{n\,}^\varepsilon \theta _{n\,}^q,\quad
\delta _n^\varepsilon \,:=\,M_n\,/\,\xi _{n\,}^\varepsilon ,\quad \delta
_{-1}^\varepsilon \,:=\,\varepsilon ,
\]
and finally 
\begin{equation}
\lambda _n^\varepsilon \;:=\;\frac{M_n}{M_{n+1}}\exp \Big[-\frac{c_0}{\theta
_n^2}\,t_n^\varepsilon \Big]  \label{notation.lambda}
\end{equation}
(with $\,c_0$\thinspace \ from Lemma \ref{lemma3.3}). We will use the
deterministic times 
\begin{equation}
T_{n+1}^\varepsilon \,:=\,T_n^\varepsilon +t_{n\,}^\varepsilon ,\quad
T_0^\varepsilon \,:=\,0.  \label{time.increments}
\end{equation}
By Remark \ref{R.decomposition} we may assume without loss of generality
that $\,\mu ({\sf R})\le 1.$\thinspace \ Then the ``starting condition'' (%
\ref{starting.cond}) is trivially satisfied.

Note that these constants satisfy Hypothesis \ref{H.stage.quantities}. In
fact, the series in condition (b2) can be estimated from above by 
\[
c\,{\rm \,}\sum_{n\ge 0}\bigg[\varepsilon \,{\rm e}^{-n}\,+\,\exp \left[
-c_{0\,}\varepsilon ^{-1}{\rm e}^{(2\alpha -\beta )n}\right] \bigg]
\]
Since $\;{\rm e}^{(2\alpha -\beta )n}\ge c\,n,$\thinspace \ this bound is of
order $\,\varepsilon ,$\thinspace \ and property (b2) follows.

By the choice of $\,\delta _{n\,}^\varepsilon ,$\thinspace \ inequality (\ref
{Feller.survival}) concerning the good paths holds by the comparison
Proposition \ref{P.comparison.Feller}.

It remains to verify the expectation estimate (\ref{small.expectation}) for
the mass of bad paths at time $\,T_{n+1\,}^\varepsilon ,$\thinspace \ for
which we will use Lemma \ref{L.new}. By time-homogeneity and definition (\ref
{E.n}) of $\,E_{n\,}^\varepsilon ,$\thinspace \ for any $\,a\in {\sf R},$%
\begin{equation}
\Pi _{T_{n\,}^\varepsilon ,a}\left( W\in _{\!_{\!_{\,}}}(E_n^\varepsilon )^{%
{\rm c}}\right) \,=\,\Pi _{0,a}\Big(\int_0^{t_n^\varepsilon }\!K({\rm d}%
r)\le \xi _n^\varepsilon \Big)  \label{Pi.expression}
\end{equation}
and 
\[
\int_0^{t_n^\varepsilon }\!K({\rm d}r)\,=\,\int_0^{t_n^\varepsilon }\!\!{\rm %
d}r\;\chi _q(W_r)\,=\,\int \!\!{\rm d}b\;\chi _q(b)\,L_{t_n^\varepsilon
}(b)\,\ge \,\int_{|b|\,\ge \,\theta _n}\!{\rm d}b\;\,\theta
_n^q\,L_{t_n^\varepsilon }(b).
\]
Thus, the probability expression in (\ref{Pi.expression}) can be estimated
from above by 
\[
\Pi _{0,a}\Big(\int_{|b|\,\ge \,\theta _n}\!{\rm d}b\;\,L_{t_n^\varepsilon
}(b)\le \frac{t_n^\varepsilon }2\Big)\;=\;\Pi _{0,a}\Big(\int_{|b|\,\le
\,\theta _n}\!{\rm d}b\;\,L_{t_n^\varepsilon }(b)\ge \frac{t_n^\varepsilon }2%
\Big).
\]
Hence, by Lemma \ref{lemma3.3}, we get 
\[
\Pi _{T_{n\,}^\varepsilon ,a}\left( W\in _{\!_{\!_{\,}}}(E_n^\varepsilon )^{%
{\rm c}}\right) \;\le \;\exp \Big[-\frac{c_0}{\theta _n^2}\,t_n^\varepsilon 
\Big]
\;=\;\frac{M_{n+1}}{M_n}\,\lambda _{n\,}^\varepsilon ,
\]
where we used (\ref{notation.lambda}). In other words, (\ref{Carl's.addition}%
) in Lemma \ref{L.new} holds, and (\ref{small.expectation}) is valid.

Altogether, we showed that all requirements for the abstract criterion
Theorem \ref{T.abstract.criterion} are satisfied, hence {\em finite time
extinction holds}\/ for the SBM with parabolic catalyst $\,\chi _q$%
\thinspace \ for any finite initial measure $\,\mu $\thinspace \ on $\,{\sf R%
}.$\hfill 
\endproof%

\subsection{A point catalyst $\protect\underline{\Gamma }$ with dense
locations\label{SS.point.dense}}

Now we consider the case 
\begin{equation}
K\,=\,L_{\left[ W,\,\underline{\Gamma }\,\right] }\quad {\rm with}\quad 
\underline{\Gamma }\;=\;\sum_{n=N}^\infty \,2^{-n\,}\underline{\pi }_{\,n\,},
\end{equation}
for a fixed $\,N\ge 0,$\thinspace \ independent of $\,\varepsilon .$%
\thinspace \ Here $\,\underline{\pi }_{\,n}$\thinspace \ is assumed to be a
(locally finite, deterministic) point measure on $\,{\sf R}$\thinspace \
such that all neighboring points have a distance of at most $\,\Delta _{n\,},
$\thinspace \ where, for some $\,\beta \in (0,1),$

\begin{equation}
\Delta _n\,:=\,{\rm e}^{-\beta n},\qquad n\ge N.  \label{Delta.n}
\end{equation}

We claim that the continuous SBM $\,X^{\underline{\Gamma }}$\thinspace \
with catalyst $\,\underline{\Gamma }$\thinspace \ {\em has the finite time
extinction property.}\/ This will follow from our abstract extinction
criterion Theorem \ref{T.abstract.criterion} once we have found the
appropriate quantities $\,T_{n\,}^\varepsilon ,M_{n\,}^\varepsilon ,\delta
_{n\,}^\varepsilon ,\lambda _{n\,}^\varepsilon ,\xi _n^\varepsilon $%
\thinspace \ entering into Hypotheses \ref{H.stage.quantities} and \ref
{H.good.and.bad}. \medskip 

\noindent $1^{\circ }$ ({\em some constants})\quad Choose $\,\alpha \in
(\frac \beta 2\,,\beta ).$\thinspace \ For $\,n\ge N$\thinspace \ and $%
\,0<\varepsilon <1,$\thinspace \ set 
\begin{equation}
m_n^\varepsilon \,:=\,\left[ \frac{{\rm e}^{\alpha n}}\varepsilon \right]
,\quad s_n^\varepsilon \,:=\,\frac{{\rm e}^{-\beta n}}{\varepsilon ^2}%
\,,\quad t_n^\varepsilon \,:=\,2\vspace{2pt}\,m_n^\varepsilon
s_{n\,}^\varepsilon ,\quad M_n^\varepsilon \,\equiv \,M_n\,:=\,2^{-n}
\label{mst}
\end{equation}
(where $\,[z]$\thinspace \ denotes the integer part of $\,z).$\thinspace \
We again use deterministic times $\,T_{n+1}^\varepsilon :=T_n^\varepsilon
+t_{n\,}^\varepsilon ,$\thinspace \ $T_N^\varepsilon :=0.$\thinspace \ Note
that by our choice of $\,t_n$\thinspace \ they satisfy (\ref
{stopping.epsilon}). The quantities $\,\xi _{n\,}^\varepsilon ,\,\delta
_{n\,}^\varepsilon ,\lambda _n^\varepsilon $\thinspace \ will be defined in (%
\ref{xi.n}), (\ref{delta.n}), and (\ref{lambda.n}), respectively.

Assume without loss of generality that $\,\mu ({\sf R})\,\le \,2^{-N}.$%
\thinspace \ Then, if we set $\,\delta _{N-1}^\varepsilon =\varepsilon ,$%
\thinspace \ the {\em starting condition}\/ (\ref{starting.cond}) is
trivially satisfied.\medskip 

\noindent 2$^{\circ }$ ({\em partitioning})\quad For $\,n\ge N$\thinspace \
and $\,0\le \varepsilon <1$\thinspace \ fixed, our next aim is to introduce
a partition of the time period $\,[T_n^\varepsilon ,T_{n+1}^\varepsilon )$%
\thinspace \ by means of some Brownian stopping times. This construction
allows us to consider hitting times of neighboring points of $\,\underline{%
\pi }_{\,n}$\ and local times spent on them.\vspace{1pt}

Given $\,\vartheta _n:=\widetilde{X}_{T_n^\varepsilon \,},$\thinspace \ and
a path $\,w$\thinspace \ ``distributed'' according to $\,\vartheta _n({\rm d}%
w),$\thinspace \ we consider the Brownian path process $\,\widetilde{W}$%
\thinspace \ distributed according to $\,\widetilde{\Pi }_{T_{n\,}^%
\varepsilon ,w\,}\,,$\thinspace \ and its projection $\,t\mapsto (\widetilde{%
W}_t)_t=W_t$\thinspace \ with law 
\[
\Pi _{T_{n\,}^\varepsilon ,w(T_n^\varepsilon )}\,=:\,\underline{\Pi }\,.
\]
(For typographical simplicity, sometimes we write $\,w(t)$\thinspace \
instead of $\,w_{t\,},$\thinspace \ etc.)\vspace{1pt}

Set $\,\kappa _0:=T_{n\,}^\varepsilon .$\thinspace \ For $\,m\ge 1,$ we
inductively define (Brownian) {\em stopping times}\/ $\,\overline{\kappa }_m=%
\overline{\kappa }_{m,n}^{\,\varepsilon }(W)$\thinspace \ and $\,\kappa
_m=\kappa _{m,n}^\varepsilon (W)$\thinspace \ as follows. Given $\,\kappa
_{m-1\,}$, let $\,\overline{\kappa }_m$\thinspace \ denote the first time
point $\,t\ge \kappa _{m-1}\,\;$such that $\,W$\thinspace \ hits one of the
atoms of $\,\underline{\pi }_{\,n\,}.$\thinspace \ Given $\,\overline{\kappa 
}_{m\,},$\thinspace \ we simply define $\,\kappa _m=\overline{\kappa }%
_m+s_n^\varepsilon $\thinspace \ [with $\,s_n^\varepsilon >0$\thinspace \
introduced in (\ref{mst})].

Write $\,H_m=H_{m,n}^\varepsilon $\thinspace \ for the {\em hitting time}\/ $%
\,\overline{\kappa }_m-\kappa _{m-1}$ of $\,\underline{\pi }_{\,n}$
(starting at time $\,\kappa _{m-1}).$\thinspace \ Recall that by definition
the distance between neighboring atoms of $\,\underline{\pi }_{\,n}$%
\thinspace \ is at most $\,\Delta _{n\,}.$\thinspace \ Using the
eigenfunction representation of solutions to the heat equation, we have that
for some constant $\,c_0>0$, 
\[
\underline{\Pi }\left( H_{m_{\!_{}}}\ge s_n^\varepsilon \right) \,\le
\,c_0^{-1}\exp \left[ -\frac{c_0\,s_n^\varepsilon }{\Delta _n^2}\right]
,\qquad m\ge 0,\quad n\ge N.
\]
Therefore, 
\begin{equation}
\underline{\Pi }\bigg(\sum_{m=1}^{m_n^\varepsilon }H_m\ge m_n^\varepsilon
s_n^\varepsilon \bigg)\le \,c_0^{-1\,}m_n^\varepsilon \exp \left[ -\frac{%
c_0\,s_n^\varepsilon }{\Delta _n^2}\right] \,=:\,\zeta _n^\varepsilon 
\label{3.3}
\end{equation}
[with $\,\Delta _{n\,},s_n^\varepsilon ,m_n^\varepsilon $\thinspace \
introduced in (\ref{Delta.n}) and (\ref{mst})].

Now write $\,L_m=L_{m,n}^\varepsilon $\thinspace \ for the (Brownian){\em \
local time}\/ spent by $\,W$\thinspace \ at the site $\,W(\overline{\kappa }%
_m)$\thinspace \ during the time interval $\,[\overline{\kappa }%
_{m\,},\kappa _m)$\thinspace \ of length $\,s_{n\,}^\varepsilon .$\thinspace
\ That is, symbolically, 
\begin{equation}
L_m:=\int_{\overline{\kappa }_m}^{\kappa _m}{\rm d}r\;\delta _{W(\overline{%
\kappa }_m)}(W_r).  \label{Y.m}
\end{equation}
Recall that at this site $\,W(\overline{\kappa }_m)$\thinspace \ there is an
atom of $\,\underline{\pi }_{\,n\,},$\thinspace \ and that the mass $\,2^{-n}
$\thinspace \ is attached to it. Therefore, using the integers $%
\,m_n^\varepsilon $\thinspace \ introduced in (\ref{mst}), for the BCLT $%
\,L_{\left[ W,\,\underline{\Gamma }\,\right] }$\thinspace \ of $\;\underline{%
\Gamma }\,\ge \,2^{-n}\,\underline{\pi }_{\,n}$\thinspace \ we get 
\begin{equation}
\int_{T_n^\varepsilon }^{\kappa _{m_n^\varepsilon }}L_{\left[ W,\,\underline{%
\Gamma }\,\right] }({\rm d}r)\,\ge \,2^{-n}\sum_{m=1}^{m_n^\varepsilon
}L_{m\,}.  \label{BCLT.bound}
\end{equation}

Clearly, the $\,L_m$\thinspace \ are i.i.d.\ (with respect to $\,\underline{%
\Pi }\,).$\thinspace \ Moreover, $\,L_m$\thinspace \ is equal in law to 
\[
\sup_{0\,\le \,t\,\le \,s_n^\varepsilon }W_t^0
\]
where $\,W^0$\thinspace \ is distributed according to $\,\Pi _{0,0}$%
\thinspace \ (see e.g.\ \cite[Theorem (6.2.3)]{RevuzYor91}). Scaling time,
we find that $\,(s_n^\varepsilon )^{-1/2}L_m$\thinspace \ is equal in law to 
$\,$%
\[
L^0:=\sup_{0\le t\le 1}W_t^0\,
\]
(which is independent of $\,n$\thinspace \ and $\,\varepsilon ).$\thinspace
\ Set $\,a:=\frac 12\,\Pi _{0,0}L^0.$\thinspace \ Since $\,L^0$\thinspace \
has finite exponential moments, by standard large deviation estimates there
exists a constant $\,c_1>0\,\,\;$such that 
\[
\underline{\Pi }\bigg((s_n^\varepsilon )^{-1/2}\sum_{m=1}^kL_m<a\,k\bigg)\le
\,{\rm e}^{-2c_1k},\qquad k\ge 1.
\]
Combining with (\ref{BCLT.bound}), we thus have 
\begin{equation}
\underline{\Pi }\bigg(\int_{T_n^\varepsilon }^{\kappa _{m_n^\varepsilon
}}L_{\left[ W,\,\underline{\Gamma }\,\right] }({\rm d}r)<\xi _n^\varepsilon %
\bigg)\le \,\exp \left[ -2c_1m_n^\varepsilon \right] ,  \label{3.2}
\end{equation}
where 
\begin{equation}
\xi _n^\varepsilon \,:=\,a\,m_n^\varepsilon (s_n^\varepsilon )^{1/2}2^{-n}\,.%
\vspace{5pt}  \label{xi.n}
\end{equation}

\noindent 3$^{\circ }$ ({\em good and bad historical paths})\quad Recall the
set $\,E_n^\varepsilon $\thinspace \ of {\em good paths}\/ introduced in
formula line (\ref{E.n}) [based on $\,t_n^\varepsilon $\thinspace \ defined
in (\ref{mst}) and entering into (\ref{time.increments}), as well as $\,\xi
_n^\varepsilon \,\;$from (\ref{xi.n})]. Since the BCLT $\,L_{\left[ W,\,%
\underline{\Gamma }\,\right] }$\thinspace \ satisfies (\ref{infinite.BCLT}),
by Proposition \ref{P.comparison.Feller} we get the survival probability
estimate (\ref{Feller.survival}) for the good paths, if we set 
\begin{equation}
\delta _n^\varepsilon \,:=\,M_n\,/\,\xi _{n\,}^\varepsilon .  \label{delta.n}
\end{equation}
On the other hand, in order to calculate the expected mass of bad paths as
required in (\ref{small.expectation}), we look at 
\begin{equation}
\underline{\Pi }\left( (E_n^\varepsilon )_{\!_{\!_{\,}}}^{{\rm c}}\right) .
\label{conditional.expect.2}
\end{equation}
In order to further estimate this, consider two cases. First let $\,w$%
\thinspace \ have ``large'' hitting times, i.e.\ 
\[
\sum_{m=1}^{m_n^\varepsilon }H_m\,\ge \,m_n^\varepsilon s_{n\,}^\varepsilon .
\]
By (\ref{3.3}), this occurs with a $\,\underline{\Pi }$--probability bounded
by $\,\zeta _{n\,}^\varepsilon .$\thinspace \ In the opposite case, by the
definition of $\,T_{n+1}^\varepsilon $\thinspace \ we have 
\[
\kappa _{m_n^\varepsilon }\,=\,T_n^\varepsilon +m_n^\varepsilon
s_n^\varepsilon +\sum_{m=1}^{m_n^\varepsilon }H_m\,<\,T_{n+1}^\varepsilon 
\]
[recall (\ref{mst}) and (\ref{time.increments})], hence here $%
\,(E_n^\varepsilon )^{{\rm c}}$\thinspace \ implies, by the definition (\ref
{E.n}) of $\,E_{n\,}^\varepsilon ,$\thinspace \ that $\,$%
\[
\int_{T_n^\varepsilon }^{\kappa _{m_n^\varepsilon }}L_{\left[ W,\,\underline{%
\Gamma }\,\right] }({\rm d}r)\;<\;\xi _{n\,}^\varepsilon .
\]
The $\,\underline{\Pi }$--probability of this event is estimated in (\ref
{3.2}). Then, for (\ref{conditional.expect.2}) we get the bound 
\begin{equation}
\zeta _{n_{\!_{\!_{}}}}^\varepsilon +\exp \left[ -2c_1m_n^\varepsilon
\right]   \label{conditional.expect.3}
\end{equation}
[with $\,\zeta _n^\varepsilon $\thinspace \ from (\ref{3.3})]. If we now set 
\begin{equation}
\lambda _n^\varepsilon \,:=\,\frac{M_n}{M_{n+1}}\left[ \zeta
_{n_{\!_{\!_{}}}}^\varepsilon +\exp \left[ -2c_1m_n^\varepsilon \right]
\right] ,  \label{lambda.n}
\end{equation}
we obtain (\ref{Carl's.addition}) in Lemma \ref{L.new} which gives (\ref
{small.expectation}).\medskip 

\noindent 4$^{\circ }$ ({\em verification of the stage quantities})\quad It
remains to check that $\,\delta _n^\varepsilon $\thinspace \ and $\,\lambda
_n^\varepsilon $\thinspace \ introduced in (\ref{delta.n}) and (\ref
{lambda.n}), respectively, satisfy Hypothesis \ref{H.stage.quantities}%
\thinspace (b2). First of all, by (\ref{delta.n}) and (\ref{xi.n}), $%
\,\delta _n^\varepsilon $\thinspace \ approximately equals $\;$%
\[
c{\rm \,}\varepsilon ^2\exp \left[ -\Big(\alpha -\frac \beta 2\Big)
n\right] ,\,\;
\]
hence its sum over $\,n$\thinspace \ is of order $\,\varepsilon ^2.$%
\thinspace \ Next, since $\,m_n^\varepsilon \ge c{\rm \,}\varepsilon ^{-1}n,$%
\thinspace \ the second term of $\,\lambda _n^\varepsilon $\thinspace \ is
bounded by $\,\exp \left[ -c{\rm \,}\varepsilon ^{-1}n\right] ,$\thinspace \
except a constant factor. Summing over $\,n$\thinspace \ we arrive at a term
of order $\,\varepsilon .$\thinspace \ Finally, by (\ref{3.3}), the first
term of $\,\lambda _n^\varepsilon $\thinspace \ is bounded from above by 
\[
\,c_0^{-1}\varepsilon ^{-1}\exp \left[ \alpha n-c_{0\,}\varepsilon ^{-2}{\rm %
e}^{\beta n}\right] .\,\;
\]
But 
\[
\varepsilon ^{-1}\int_1^\infty \!{\rm d}x\;\exp \left[ \alpha
x-c_{0\,}\varepsilon ^{-2}{\rm e}^{\beta x}\right] \;\le \;\left( \beta
\varepsilon \right) ^{-1}\int_1^\infty \!{\rm d}y\;y^{\frac \alpha \beta
-1}\exp \left[ -c_{0\,}\varepsilon ^{-2}y\right] ,
\]
and $\,r\mapsto r^z{\rm e}^{-r}$\thinspace \ is bounded for $\,r>0$%
\thinspace \ bounded away from $\,0,$\thinspace \ for each fixed $\,z.$%
\thinspace \ Hence, it suffices to consider 
\[
\varepsilon ^{-1}\int_1^\infty \!{\rm d}y\;\exp \left[ -\frac{c_0}%
2\varepsilon ^{-2}y\right] 
\]
which is of order $\,\varepsilon .$\thinspace \ Consequently, the series in
Hypothesis \ref{H.stage.quantities}\thinspace (b2) is bounded by $\;c{\rm \,}%
\varepsilon ,$\thinspace \ hence this hypothesis is satisfied in the present
case.

Summarizing, the {\em catalytic SBM }$\,X^{\underline{\Gamma }}$\thinspace \ 
{\em dies in finite time,}\/ for any finite starting measure $\,\mu $%
\thinspace \ on $\,{\sf R}.$\hfill 
\endproof%

\section{Proof of the main result\label{S.proof.main}}

The proof of Theorem \ref{T.finite.time.exti} (p.\pageref{T.finite.time.exti}%
) concerning finite time extinction of the one-dimensional SBM $\,X^\Gamma $%
\thinspace \ with a stable catalyst proceeds in several steps. Since here we
start with an initial measure $\,\mu $\thinspace \ of compact support, and $%
\,X^\Gamma \,\,\;$has the compact support property (\cite{DawsonLiMueller95}%
), by some coupling technique we will pass to a periodic catalyst $\,\Gamma
^{{\sf K}}$. Moreover, because the survival probability is monotone in the
catalyst, we will switch to a smaller catalyst, as already explained in
Subsection \ref{SS.Model1}. Altogether we will reduce to the case of a point
catalyst $\,\underline{\Gamma }$\thinspace \ with dense locations as dealt
with in Subsection \ref{SS.point.dense}.

\subsection{A coupling of catalytic SBMs}

Recall that the historical catalytic SBM $\,\widetilde{X}^\Gamma $\thinspace
\ exists for $\,%
\bP%
$--almost all $\,\Gamma .$\thinspace \ Fix an initial measure $\,\mu \in 
{\cal M}_{{\rm f}}$\thinspace \ with compact support. We want to show that 
\begin{equation}
\widetilde{P}_{0,\mu }^\Gamma \left( \widetilde{X}_t^\Gamma \ne 0,\;\forall
t\right) \;=\;0,\quad 
\bP%
{\rm -a.s.}  \label{claim.tilde}
\end{equation}

For $\,{\sf K}\ge 1,$\thinspace \ let 
\[
E_{{\sf K}}\;:=\;\left\{ w\in {\bf C}:_{\!_{\!_{\,_{}}}}\;|w_s|\le {\sf K}%
,\;\forall s\ge 0\right\} .
\]
According to the compact support property of \cite{DawsonLiMueller95}, 
\begin{equation}
\lim_{{\sf K}\uparrow \infty }\,\widetilde{P}_{0,\mu }^\Gamma \left( {\rm %
supp}\widetilde{X}_t^\Gamma \subseteq E_{{\sf K\,}},\;\forall t\right)
\;=\;1,\quad 
\bP%
{\rm -a.s.}  \label{csp}
\end{equation}
For the further proof, fix such a sample $\,\Gamma .$\thinspace \ By (\ref
{csp}), instead of (\ref{claim.tilde}) it suffices to show that 
\begin{equation}
\widetilde{P}_{0,\mu }^\Gamma \left( \widetilde{X}_t^\Gamma \ne 0\;{\rm %
and\;\,supp}\widetilde{X}_t^\Gamma \subseteq E_{{\sf K\,}},\;\forall
t\right) \;=\;0,\quad {\rm for\;all\;\,}{\sf K}.  \label{claim.K}
\end{equation}
But under this restriction to historical paths living in $\,E_{{\sf K\,}},$%
\thinspace \ we may change the catalyst outside of $\,[-{\sf K},{\sf K}]$%
\thinspace \ without affecting the latter probability. This will be
formalized in the following considerations establishing some coupling
argument.

Fix $\,{\sf K}\ge 1$\thinspace \ such that the initial measure $\,\mu $%
\thinspace \ is supported by $\,(-{\sf K},{\sf K}).$\thinspace \ Consider
the hitting time $\,\tau ^{{\sf K}}$\thinspace \ of $\,\left\{ -{\sf K},{\sf %
K}\right\} ,$\thinspace \ the boundary of the interval $\,\left( -{\sf K},%
{\sf K}\right) .$ Then replace the Brownian motion \thinspace $W$\thinspace
\ of reactant particles by the stopped process $\,t\mapsto W_{t\wedge \tau ^{%
{\sf K}}\,}.$\thinspace \ This \vspace{2pt}transfers $\,\widetilde{X}^\Gamma 
$\thinspace \ into the stopped historical catalytic SBM $\,\vspace{1pt}%
t\mapsto \widetilde{X}_{t\wedge \tau ^{{\sf K}}\,}^\Gamma .$\thinspace \
Note that the paths of this stopped process live completely in the closed
interval $\,\left[ -{\sf K},{\sf K}\right] .$\thinspace \ 

Actually we {\em decompose}\/ this stopped process, 
\[
\widetilde{X}_{t\wedge \tau ^{{\sf K}}}^\Gamma \;=\;{\rm m}_t^{{\sf K}}+{\rm %
m}_{t\,}^{\circ },
\]
by distinguishing between paths $\,$%
\[
w\,\in \,^{{\sf K}\!}{\bf C}^t\,:=\,\left\{ w\in {\bf C}^t:\;\left|
w_t\right| ={\sf K}\right\} \ 
\]
which end at the boundary $\,\left\{ -{\sf K},{\sf K}\right\} ,$ \thinspace
\ and those which stay within $\,\left( -{\sf K},{\sf K}\right) :$%
\[
w\,\in \,^{\circ \!}{\bf C}^t\,:=\,\left\{ w\in {\bf C}^t:_{\!_{\!_{\,_{}}}}%
\;\left| w_s\right| <{\sf K},\;\forall s\le t\right\} .
\]
In other words, 
\[
{\rm m}_t^{{\sf K}}\;:=\;\widetilde{X}_{t\wedge \tau ^{{\sf K}}}^\Gamma
\left( \left( \cdot \right) \,_{\!_{\!_{\,}}}\cap \,^{{\sf K}\!}{\bf C}%
^t\right) ,\qquad {\rm m}_t^{\circ }\;:=\widetilde{X}_{t\wedge \tau ^{{\sf K}%
}}^\Gamma \left( \left( \cdot \right) \,_{\!_{\!_{\,}}}\cap \,^{\circ \!}%
{\bf C}^t\right) .
\]
Note that the path $\,t\mapsto {\rm m}_t^{{\sf K}}$\thinspace \ of measures
on $\,{\bf C}$\thinspace \ is monotonically non-decreasing. Thus, $\,{\rm m}%
_{{\rm d}t}^{{\sf K}}({\rm d}\omega _t)$\thinspace \ can be considered as a
measure in $\,{\cal M}\!\left[ {\sf R}_{+}\widehat{\times }{\bf C}^{\bullet
}\right] $\thinspace \ [recall notation (\ref{cross.hut})]. Now we use the
increments of this historical path $\,{\rm m}^{{\sf K}}$\thinspace \ as an 
{\em immigration}\/ process of a historical catalytic SBM starting from the
zero measure, denoted by $\,\widetilde{Y}=\widetilde{Y}^{\Gamma ,{\rm m}^{%
{\sf K}}}.\,\,$More precisely, (for the given $\,\Gamma )$\thinspace \
given\/ $\,\,{\rm m}^{{\sf K}},$\thinspace \ defining $\,\widetilde{Y}=%
\widetilde{Y}^{\Gamma ,{\rm m}^{{\sf K}}}$\thinspace \ we use the modified
process according to Proposition \ref{P.cont.historical.SBM}\thinspace (b)
with the collision local time $L_{[W,\Gamma ]}$\thinspace \ as branching
rate functional, and with $\,\eta $\thinspace \ defined by 
\[
\eta ({\rm d}r,{\rm d}\omega _r)\;:=\;{\rm m}_{{\rm d}r}^{{\sf K}}({\rm d}%
\omega _r).
\]

We need also another process. Let $\,\Gamma ^{{\sf K}}$\thinspace \ denote
the {\em periodic extension}\/\ of the restriction $\,\Gamma \left( \left(
\cdot \right) \cap (-{\sf K},{\sf K}]\right) $\thinspace \ of $\,\Gamma $%
\thinspace \ to $\,(-{\sf K},{\sf K}]$\thinspace \ to all of $\,{\sf R}$%
\thinspace \ (for the fixed $\,\Gamma ).$\thinspace \ Now replace $\,\Gamma $%
\thinspace \ by $\,\Gamma ^{{\sf K}}$\thinspace \ in the definition of $\,%
\widetilde{Y}=\widetilde{Y}^{\Gamma ,{\rm m}^{{\sf K}}},$\thinspace \ to
obtain a historical catalytic SBM with {\em periodic}\/ catalyst $\,\Gamma ^{%
{\sf K}}$\thinspace \ and {\em immigration}\/ controlled by $\,{\rm m}^{{\sf %
K}},$\thinspace \ which we denote by $\,\widetilde{Z}=\widetilde{Z}^{\Gamma
^{{\sf K}},{\rm m}^{{\sf K}}}.$

Recall that both of our processes $\,\widetilde{Y}$\thinspace \ and $\,%
\widetilde{Z}$\thinspace \ are based on the same samples $\,\,\Gamma \left(
\left( \cdot \right) \cap (-{\sf K},{\sf K}]\right) $\thinspace \ and $\,%
{\rm m}^{{\sf K}}.$\thinspace \ The reason we introduced these processes is
the following obvious coupling result.

\begin{lemma}[coupling of historical catalytic SBMs]
\label{L.coupling}Fix\/ $\,{\sf K}\ge 1$\thinspace \ such that the initial
measure $\,\mu $\thinspace \ is supported by $\,(-{\sf K},{\sf K}).$%
\thinspace \ Given\/ $\,\,\Gamma ,$\thinspace \ the processes $\,\,%
\widetilde{Y}+{\rm m}^{\circ }$\thinspace \ and $\,\widetilde{Z}+{\rm m}%
^{\circ }$\thinspace \ coincide in law with the historical catalytic SBMs $\,%
\widetilde{X}^\Gamma $\thinspace \ and $\,\widetilde{X}^{\Gamma ^{{\sf K}}}$%
\thinspace \ with catalysts $\,\Gamma $\thinspace \ and $\,\Gamma ^{{\sf K}},
$\thinspace \ respectively, and their restrictions to paths living in $\,E_{%
{\sf K}}$\thinspace \ are identical to \thinspace ${\rm m}^{\circ }$.
\end{lemma}

\subsection{Completion of the proof of the main theorem}

By Lemma \ref{L.coupling}, we may pass in (\ref{claim.K}) from $\,\Gamma $%
\thinspace \ to the periodic $\,\Gamma ^{{\sf K}}.$\thinspace \ Hence,
instead of (\ref{claim.K}) it suffices to show 
\begin{equation}
\widetilde{P}_{0,\mu }^{\Gamma ^{{\sf K}}}\left( \widetilde{X}_t^{\Gamma ^{%
{\sf K}}}\ne 0,\;\forall t\right) \;=\;0,  \label{claim.periodic}
\end{equation}
for each fixed $\,{\sf K}\ge 1$\thinspace \ such that $\,\mu $\thinspace \
is supported by $\,(-{\sf K},{\sf K}).$\thinspace \ In other words, we want
to show finite time extinction of the historical catalytic SBM $\,\widetilde{%
X}^{\Gamma ^{{\sf K}}}$\thinspace \ with fixed periodic catalyst $\,\Gamma ^{%
{\sf K}}.$

In order to can apply later on the result of Subsection \ref{SS.point.dense}%
, we further use the fact that the collision local times $\,L_{[W,\psi ]}$%
\thinspace \ are non-decreasing in $\,\psi \in {\cal M}_{p\,}.$\thinspace \
That is, $\,\psi _1\le \psi _2$\thinspace \ implies $\,L_{[W,\psi _1]}\le
L_{[W,\psi _2]}.$\thinspace \ Therefore the corresponding solutions $%
\,v^\psi $\thinspace \ of the log-Laplace equation (\ref{log.Laplace.X.hist}%
) are non-increasing: $\,v^{\psi _1}\ge v^{\psi _2}.$\thinspace \ But this
yields that the {\em extinction probability is non-decreasing}\/ in $\,\psi :
$%
\[
\psi _1\,\le \,\psi _2\quad {\rm implies}\quad P_{0,\mu }^{\psi
_1}(X_t=0)\,\le \,P_{0,\mu }^{\psi _2}(X_t=0)
\]
[recall (\ref{exti.prob})]. Hence, for our purpose of verifying (\ref
{claim.periodic}), we may replace the periodic catalyst $\,\Gamma ^{{\sf K}}$%
\thinspace \ by a smaller measure.

To this end, as already mentioned in Subsection \ref{SS.Model1}, we first
drop all the ``big'' point catalysts: For the moment, fix $\,N\ge 0$%
\thinspace \ (independent of $\,\varepsilon ),$\thinspace \ and remove all
those atoms $\,\alpha _{i\,}\delta _{b_i}$\thinspace \ of $\,\Gamma ^{{\sf K}%
}$\thinspace \ [or $\,\Gamma ,$\thinspace \ recall the representation (\ref
{Gamma.rep})] with action weight $\,\vspace{1pt}\alpha _i\ge 2^{-N+1}.$%
\thinspace \ Next, for each $\,n\ge N,\,\,\;$we replace the action weights $%
\,\alpha _i\in [2^{-n},2^{-n+1})\,\;$by $\,2^{-n}.\,\,\;$Note that with
respect to $\,%
\bP%
,$\thinspace \ the positions $\,b_i\in (-{\sf K},{\sf K}]$\thinspace \ of
the related atoms are distributed as a {\em Poisson point process}\/ with
intensity measure $\,c_\gamma \,2^{\gamma n}\,{\bf 1}_{(-{\sf K},{\sf K}%
]}(b)\,\ell ({\rm d}b).$\thinspace \ Here the constant $\,c_\gamma $%
\thinspace \ is given by 
\[
c_\gamma \,:=\,\,\gamma ^{-1}(1-2^{-\gamma })\Big(\int_0^\infty {\rm d}%
r\;\,r^{-1-\gamma }\,(1-{\rm e}^{-r})\Big)^{-1}
\]
(see e.g.\ \cite{DawsonFleischmann92.equ}). Let \vspace{1pt}$\,\underline{%
\pi }_{\,n}$\thinspace \ denote the periodic extension of this Poisson point
process, extension from $\,(-{\sf K},{\sf K}]$\thinspace \ to all of $\,{\sf %
R}.$

What remains for the reduction to Subsection \ref{SS.point.dense} is to show
that $\,%
\bP%
$--a.s.\ in $\,\underline{\pi }_{\,n}$\thinspace \ neighboring catalysts
have a distance of at most $\,\Delta _n={\rm e}^{-\beta n},$\thinspace \ for
all $\;n\ge N,$\thinspace \ for $\,N\ge 0$\thinspace \ appropriately chosen.
For this purpose, we fix 
\begin{equation}
\beta \in \left( 0,\gamma \log 2\right) .  \label{beta}
\end{equation}
By Borel-Cantelli, it suffices to show that the quantities 
\begin{equation}
\bP%
\left( \exists \;{\rm two\;neighboring\;points\;in\;\,}\underline{\pi }%
_{\,n}\,\;{\rm with\;a\;distance\;larger\;than\;}\,_{\!_{\!_{\,_{}}}}{\rm e}%
^{-\beta n}\right)   \label{summable}
\end{equation}
are summable in $\,n\ge 1.$\thinspace \ But each of these probabilities is
bounded from above by 
\begin{equation}
\bP%
\left( \max_{1\,\le \,i\,\le \,J_n+1}\xi _i\,>{\rm e}^{-\beta n}\,\right) 
\label{summable'}
\end{equation}
where $\,J_n$\thinspace \ is the Poissonian number of points in $\,(-{\sf K},%
{\sf K}]$\thinspace \ with expectation $\,a_n:=2{\sf K}c_\gamma 2^{\gamma n},
$\thinspace \ and the $\,\xi _{1\,},\xi _{2\,},...$\thinspace \ are i.i.d.\
exponentials with parameter $\,a_{n\,}.$\thinspace \ Now the $\,J_n$%
\thinspace \ satisfy a standard large deviation principle as $\,n\uparrow
\infty ,$\thinspace \ hence, 
\[
\bP%
\left( J_n+1_{\!_{\!_{\,}}}>2a_n\right) \;\le \;\exp \left[ -\,c{\rm \,}%
2^{\gamma n}\right] 
\]
for all sufficiently large $\,n.$\thinspace \ Since the right hand side is
summable in $\,n,$\thinspace \ in the probability expression (\ref{summable'}%
) we may additionally restrict to $\,J_n+1\le 2a_{n\,}.$\thinspace \
Consequently, instead of (\ref{summable'}) we look at 
\[
\bP%
\left( \max_{1\,\le \,i\,\le \,2a_n}\xi _i\,>{\rm e}^{-\beta n}\,\right) .
\]
By scaling, we may switch to 
\begin{equation}
{\cal P}\left( \max_{1\,\le \,i\,\le \,2a_n}\xi _i^{\prime }\,>a_n\,{\rm e}%
^{-\beta n}\,\right)   \label{summable''}
\end{equation}
where the $\,\xi _i^{\prime }$\thinspace \ are now i.i.d.\ standard
exponentials (under the law denoted by $\,{\cal P}).$

Next we use the fact that for all $\,x\ge 0$\thinspace \ and $\,m\ge 2,$%
\begin{equation}
\left| {\cal P}\left( \max_{1\,\le \,i\,\le \,m}\xi _i^{\prime }-\log
m\,>x\right) -\left( 1-_{\!_{\!_{\,}}}\exp \left[ -{\rm e}^{-x}\right]
\right) \right| \;\le \;2\,{\rm e}^{-2x}
\end{equation}
(see Example 2.10.1 in \cite{Galambos1978}; take $\,q=\frac 12$\thinspace \
there). Now, for all $\,n$\thinspace \ sufficiently large, $\,m=\left[
2a_n\right] $\thinspace \ and 
\begin{equation}
x\;=\;a_n\,{\rm e}^{-\beta n}-\log \left( 2a_n\right) \;=\;c{\rm \,e}%
^{n\left( \gamma \log 2-\beta \right) }-\log (4{\sf K}c_\gamma )-n\gamma
\log 2  \label{x}
\end{equation}
satisfy these conditions [recall (\ref{beta})]. Thus, for (\ref{summable''})
we get the bound 
\[
1-\exp \left[ -{\rm e}^{-x}\right] +2\,{\rm e}^{-2x}\;\le \;3\,{\rm e}^{-x}
\]
which for $\,x$\thinspace \ from (\ref{x}) is summable in $\,n.$

This finishes the proof of Theorem \ref{T.finite.time.exti}.\hfill 
\endproof%

\section{The lattice model\label{SS.lattice}}

Now consider the model with random catalysts on the lattice $\,{\sf Z}^d.$%
\thinspace \ Recall that $\,\varrho =\{\varrho _b\}_{b\in {\sf Z}^d}$%
,\thinspace $\;$the catalysts, are i.i.d.\/\ random variables which are
uniformly distributed on $\,[0,1].$\thinspace \ Instead of Brownian motions,
the motion process is now given by a continuous time simple random walk on $%
\,{\sf Z}^d,$\thinspace \ which moves to a neighboring site at rate 1. In
other words, the times between jumps are i.i.d.\ exponential times with mean
1.

We use symbols analogous to the ones in earlier sections. In particular, $\,%
\bP%
$\thinspace \ denotes the law of the catalyst, $W=\left[ W,\,\Pi
_{s,a\,},\,s\ge 0,\,\;a\in {\sf Z}^d\right] $\thinspace \ the {\em simple
random walk}\/ in $\,{\sf Z}^d$\thinspace \ on canonical Skorohod path space 
$\,{\bf D}={\cal D}\left[ {\sf R}_{+\,},{\sf Z}^d\right] $\thinspace \ of
c\`{a}dl\`{a}g functions, and $\,$%
\[
\widetilde{X}^\varrho =\left[ \widetilde{X}^\varrho ,\,\widetilde{P}_{s,\mu
\,}^\varrho ,\,s\ge 0,\,\mu \in {\cal M}_{{\rm f}}[{\bf D}^s]\right] 
\]
the {\em historical simple super-random walk}\/ on $\,{\sf Z}^d$\thinspace \
(also called simple interacting Feller's branching diffusion) with catalyst $%
\,\varrho .$\thinspace \ Note that Proposition \ref{P.cont.historical.SBM},
Theorem \ref{T.abstract.criterion}, Lemma \ref{L.new}, and Proposition \ref
{P.comparison.Feller} remain valid (with the obvious changes).

For simplicity, we now assume that $\,X_0=\delta _{0\,}.$\thinspace \ Our
aim is to show the {\em finite time extinction}\/ property for $\,\widetilde{%
X}^\varrho ,$\thinspace \ for $\,%
\bP%
$--a.a.\ $\,\varrho .$\thinspace \ In this case, the bad historical paths
are those which spend a large amount of time at sites $\,b\in {\sf Z}^d$%
\thinspace \ where $\,\varrho _b$\thinspace \ is small. We will choose time $%
\,T_N$\thinspace \ so that, with high probability, most of the mass is dead
by this time [in the sense of (\ref{starting.cond})]. This is the hardest
part of the argument. Bounding the mass after this uses similar but easier
ideas.

We also need the following crude estimate on the {\em distance traveled by
the simple random walk}\/ $\,W$\thinspace \ in time $\,t.$\thinspace \ Let $%
\,J_t$\thinspace \ denote the number of jumps taken by $\,W$\thinspace \ by
time $\,t.$\thinspace \ Since $\,J_t$\thinspace \ is {\em Poisson}\/ with
parameter $\,t,$\thinspace \ for $\,k\ge 0$\thinspace \ we have 
\begin{equation}
\left. 
\begin{array}{l}
\Pi _{0,0}\left( \sup\limits_{0\le s\le t}\left| W_s\right| \ge k\right)
\;\le \;\Pi _{0,0}\left( J_t\ge k\right) \;=\;\,{\rm e}^{-1}\dsum%
\limits_{i=k}^\infty \dfrac{t^i}{i!}\vspace{4pt} \\ 
\;\le \;{\rm e}^{-1}\,\dfrac{t^k}{k!}\,\dsum\limits_{i=0}^\infty \dfrac{t^i}{%
i!}\;=\;\dfrac{t^k}{k!}\;\le \;\left( \dfrac{t{\rm e}}k\right) ^{\!k}(2\pi
k)^{-1/2},
\end{array}
\;\right\}   \label{3.14}
\end{equation}
the latter by Stirling's approximation.

For $\,n\ge 0,$\thinspace \ let $\,D_n$\thinspace \ denote the cube 
\begin{equation}
D_n\,=\,\left\{ (b_1,\dots ,b_d)\in {\sf Z}^d:_{\,_{\,_{\,_{}}}}\max
(|b_1|,\dots ,|b_d|)\le 2^n\right\} 
\end{equation}
in $\,{\sf Z}^d$\thinspace \ having $\,(2^{n+1}+1)^d$\thinspace \ sites. For
a given path $\,w\in {\bf D},$\thinspace \ let $\,\tau _n=\tau _n(w)$%
\thinspace \ denote the first time $\,t\ge 0$\thinspace \ that $\,w_t$%
\thinspace \ does not belong to $\,D_{n\,}.$\thinspace \ We intend to use 
{\em Dynkin's special Markov property}\/ to start the stopped historical
super-random walk $\,\left\{ \widetilde{X}_{\tau _n}^\varrho :\;n\ge
0\right\} $\thinspace \ afresh at the times $\,\tau _{n\,}.$

We define in this proof that the quantities $\,M_{n\,}^\varepsilon
,\,\lambda _{n\,}^\varepsilon ,\,\xi _{n\,}^\varepsilon ,\,\delta
_{n\,}^\varepsilon ,\,T_{n\,}^\varepsilon ,$\thinspace $\;$and $%
\,E_n^\varepsilon $\thinspace \ entering in Hypotheses \ref
{H.stage.quantities} and \ref{H.good.and.bad} to be independent of $%
\,\varepsilon ,$\thinspace \ and therefore we omit the index $\varepsilon
.\,\;$We will choose $\,N=N(\varepsilon )\ge 1$\thinspace \ later, such that 
$\,\lim_{\varepsilon \downarrow 0}N(\varepsilon )=\infty .$\thinspace \ To
be more specific, $\,N$\thinspace \ must be so large that all of the
statements involving the phrase ``for $\,N$\thinspace \ sufficiently large''
are satisfied. Set
\[
\delta _{N-1}\,:=\,2^{-N/4},\quad T_N\,:=\,\frac{2^N}6\wedge \tau _{N\,},
\]
and for $\,n\ge N=N(\varepsilon ),$\thinspace \ let 
\[
M_n\,:=\,2^{-n(d+3)},\quad \lambda _n\,:=\,2^{-2^n},\quad \delta
_n\,:=\,2^{-n},
\]
\[
\xi _n\,:=\,2^{-n(d+2)},\quad T_{n+1}\,:=\,\left( T_n+2^{-n}\right) \wedge
\tau _{n+1\,}.
\]
Note that $\,\delta _{N-1}$\thinspace \ and $\,T_N$\thinspace \ implicitly
depend on $\,\varepsilon $\thinspace \ via $\,N(\varepsilon ).$\thinspace \
One can easily show that $\,T_{n+1}>T_{n\,},$\thinspace \ $\Pi _{0,0}$--a.s.
Clearly (b1) and (b2) of Hypothesis \ref{H.stage.quantities} are satisfied.

We need the following large deviations lemma on the simple random walk $\,W.$%
\thinspace \ We say that a non-empty subset $\,S\subset {\sf Z}^d$\thinspace
\ is {\em connected}\/ if any two elements $\,a,b\in S$\thinspace \ are
connected by a chain $\,a=z_0,\dots ,z_k=b$\thinspace \ of elements of $\,S,$%
\thinspace \ such that for $\,1\le i\le k$\thinspace \ the points $%
\,z_{i-1},z_i$\thinspace \ are nearest neighbors. That is, they are distance 
$1$ apart.

\begin{lemma}[large deviations]
\label{lemma.clusters}Fix\/ $\,m\ge 1.\,\,\;$Suppose that\/ $\,S\subset {\sf %
Z}^d$\thinspace \ has the property that no connected subset of\/ $\,S$%
\thinspace \ has cardinality larger than\/ $\,m{.\,\,\;}$Then there exist
constants\/ $\,\alpha ,c>0$\thinspace \ {\em (}depending on\/ $\,m{)}$%
\thinspace \ such that for all\/ $\,t\ge 1$, 
\[
\sup_{a\in {\sf Z}^d}\,\Pi _{0,a}\left( \int_0^t\!{\rm d}s\;{\bf 1}\left\{
W_s\in S^{{\rm c}}\right\} \le \alpha t\right) \,\le \;c^{-1}{\rm e}^{-ct}.
\]
\end{lemma}

\proof%
By monotonicity in $\,m,$\thinspace \ we may enlarge $\,m$\thinspace \ if
necessary, so we may assume that $\,m$\thinspace \ is even. By our
assumptions, if $\,a\in S,$\thinspace \ then there exists a chain consisting
of points $\,a=z_0,\ldots ,z_m$\thinspace \ such that $\,z_{i-1\,},z_i$%
\thinspace \ are nearest neighbors for $\,1\le i\le m$, and $\,z_m=z_m(a)\in
S^{{\rm c}}.$\thinspace \ If $\,a\in S^{{\rm c}},$\thinspace \ we construct
such a chain as follows. Let $\,b$\thinspace \ be one of the nearest
neighbors of $\,a,$\thinspace \ and let $\,z_{2k}:=a,$\thinspace \ $%
z_{2k+1}:=b.$\thinspace \ 

Suppose that $\,W_0=a.\,\,\;$Let $\,\eta _a$\thinspace \ denote the first
time $\,t$\thinspace \ that $\,W_t=z_m(a).$\thinspace \ (If there is no such
time, set $\,\eta _a=\infty .)$\thinspace \ Let $\,F=F(a)$\thinspace \
denote the event that $\,\eta _a<1/2$\thinspace \ and that $\,W_s=z_m(a)$%
\thinspace \ for $\,\eta _a\le s\le 1.$\thinspace \ By the properties of our
continuous-time simple random walk, using the constructed chain, there
exists $\,\alpha >0$\thinspace \ such that for all $\,a\in {\sf Z}^d,$%
\begin{equation}
\Pi _{0,a}\left( F(a)_{\!_{\!_{\,}}}\right) \,\ge \;8\alpha .  \label{varep}
\end{equation}
Let $\,F_i:=\theta _iF,$\thinspace \ $i\ge 0,$\thinspace \ where $\,\theta _s
$\thinspace \ is the time-shift operator on paths. By the Markov property
and (\ref{varep}), there exists a sequence of independent events $\,%
\overline{F}_i$\thinspace \ such that $\,\overline{F}_i\subset F_i$%
\thinspace \ and\ $\,\Pi _{0,a}(\overline{F}_i)=8\alpha ,$\thinspace \ for
each $\,i.$\thinspace \ Set $\,$%
\[
G_k\,:=\,\sum_{i=0}^{k-1}{\bf 1}_{\overline{F}_i\,}.
\]
By Chernoff's large deviations theorem (see \cite[Theorem 9.3]{Billingsley86}%
), there exist a constant $\,c>0$\thinspace \ such that for all $\,a\in {\sf %
Z}^d,$\thinspace \ we have 
\[
\Pi _{0,a}\left( \frac{G_k}k\le 4\alpha \right) \,\le \,c^{-1}{\rm e}%
^{-ck},\qquad k\ge 0.
\]
Note that 
\begin{equation}
\int_i^{i+1}\!{\rm d}s\;{\bf 1}\left\{ W_s\in S^{{\rm c}}\right\} \,\ge
\,1/2\quad {\rm on\;\,}F_{i\,}.  \label{display}
\end{equation}
Indeed, if $\,F_i$\thinspace \ occurs, then $\,W_s\in S^{{\rm c}}$\thinspace
\ for $\,s\in \left( i+1/2\,,i+1\right) .$\thinspace \ Suppose that $%
\,G_k/k\ge 4\alpha .$\thinspace \ Then there are at least $\,4\alpha k$%
\thinspace \ indices $\,i\le k-1$\thinspace \ such that $\,\overline{F}_i$%
\thinspace \ occurs, and hence $\,F_i$\thinspace \ occurs. In that case, by (%
\ref{display}), 
\[
\int_0^k\!{\rm d}s\;{\bf 1}\left\{ W_s\in S^{{\rm c}}\right\} \;\ge
\;2\alpha k.
\]
Hence, for $\,a\in {\sf Z}^d$\thinspace \ and $\,k\ge 0,$%
\[
\Pi _{0,a}\left( \frac{G_k}k>4\alpha \right) \;\le \;\Pi _{0,a}\left(
\int_0^k\!{\rm d}s\;{\bf 1}\left\{ W_s\in S^{{\rm c}}\right\} \;\ge
\;2\alpha k\right) .
\]
Interpolating, we have that for all $\,a\in {\sf Z}^d$\thinspace \ and $%
\,t\ge 1$, 
\begin{eqnarray*}
\Pi _{0,a}\left( \int_0^t\!{\rm d}s\;{\bf 1}\left\{ W_s\in S^{{\rm c}%
}\right\} <\alpha t\right) \!\! &\le &\!\!\Pi _{0,a}\left( \int_0^{[t]}\!%
{\rm d}s\;{\bf 1}\left\{ W_s\in S^{{\rm c}}\right\} <\alpha \left(
[t]+1_{\!_{\!_{\,}}}\right) \right)  \\
\!\! &\le &\!\!c^{-1}\,{\rm e}^{-ct},
\end{eqnarray*}
finishing the proof of Lemma \ref{lemma.clusters}.\hfill 
\endproof%
\bigskip 

For $\,m\ge 1,$\thinspace \ $n\ge N=N(\varepsilon ),$\thinspace \ and $%
\,0\le \zeta \le 1,$\thinspace \ let $\,A(m,n,\zeta )$\thinspace \ denote
the (catalyst) event that there is no connected subset $\,S\subset D_n$%
\thinspace \ with cardinality $\,m,$\thinspace \ on which all of the
catalysts satisfy $\,\varrho _b\le \zeta .$\thinspace \ Note that there is a
finite number $\,c(m,d)$\thinspace \ of connected sets of cardinality $%
\,m,\,\,\;$which contain a given point. Then we have 
\begin{equation}
\left. 
\begin{array}{lll}
\bP%
\left( A^{{\rm c}}(m,n,\zeta )_{\!_{\!_{\,}}}\right)  & \le  & \left(
2^{n+1}+1\right) ^dc(m,d)\,\left( 
\bP%
\left( \varrho _b\le \zeta \right) _{\!_{\!_{\,}}}\right) ^m\vspace{8pt} \\ 
& = & \left( 2^{n+1}+1\right) ^dc(m,d)\,\zeta ^m.
\end{array}
\;\right\}   \label{clusters}
\end{equation}
In particular, if $\;$%
\begin{equation}
\zeta =\zeta _n=2^{-(n-1)(d+1)},  \label{notation.zeta}
\end{equation}
then 
\begin{equation}
\bP%
\left( A^{{\rm c}}(1,n,\zeta _n)_{\!_{\!_{\,}}}\right) \,\le
\,c\,2^{nd\,}2^{-(n-1)(d+1)}\,=\,c\,2^{-n}.  \label{cluster2}
\end{equation}

For $\,m=1,$\thinspace \ all catalysts in $\,D_n$\thinspace \ are greater
than $\,\,\zeta $\thinspace \ on $\,A(1,n,\zeta ).$\thinspace \ Put 
\begin{equation}
A_1(n)\,:=\,\bigcap_{k=n}^\infty A(1,k,\zeta _k)
\end{equation}
and note that 
\begin{equation}
\bP%
\left( A_1^{{\rm c}}(n)_{\!_{\!_{\,}}}\right) \,\le \,c_1\,2^{-n}.
\label{clusters3}
\end{equation}

From now on, let $\,$%
\begin{equation}
m=2(d+1),  \label{notation.m}
\end{equation}
and take $\,$%
\begin{equation}
\overline{\zeta }=\overline{\zeta }_n=2^{-n/2}.  \label{notation.zeta.bar}
\end{equation}
Then, by (\ref{clusters}), we have, 
\[
\bP%
\left( A^{{\rm c}}(m,n,\overline{\zeta }_n)_{\!_{\!_{\,_{\!}}}}\right)
\,<\,c\,2^{-n}.
\]
Let 
\begin{equation}
A_2(n)\,:=\,\bigcap_{k=n}^\infty A\left( m,k,\overline{\zeta }_k\right) 
\end{equation}
and note that 
\begin{equation}
\bP%
\left( A_2^{{\rm c}}(n)_{\!_{\!_{\,}}}\right) \,\le \,c_2\,2^{-n}.
\label{clusters4}
\end{equation}
Fix $\,\overline{\varepsilon }>0.$\thinspace \ Using (\ref{clusters3}) and (%
\ref{clusters4}), we choose $\,\overline{n}=\overline{n}(\overline{%
\varepsilon })$\thinspace \ so large that 
\begin{equation}
\bP%
\left( A_1^{{\rm c}}(\overline{n})\cup A_{2_{\!_{}}}^{{\rm c}}(\overline{n}%
)\right) \;\le \;\left( c_1+c_2\right) 2^{-\overline{n}}\;<\;\overline{%
\varepsilon }.
\end{equation}
We will apply our general Theorem \ref{T.abstract.criterion} with $\,N$%
\thinspace \ chosen to satisfy $\,N=N(\varepsilon ,\overline{\varepsilon }%
)\ge \overline{n}.$\thinspace \ We will conclude that for catalysts $%
\,\varrho $\thinspace \ in $\,\,A_1(\overline{n})\cap A_2(\overline{n}),$%
\thinspace \ finite time extinction occurs with $\,\widetilde{P}_{0,\mu
}^\varrho $--probability 1. Therefore, with $\,%
\bP%
$--probability at least $\,1-\overline{\varepsilon },$\thinspace \ finite
time extinction occurs. Since $\,\overline{\varepsilon }$\thinspace \ is
arbitrary, our proof will then be finished.\smallskip 

From now on we assume that $\,\varrho $\thinspace \ belongs to the set $%
\,A_1(\overline{n})\cap A_2(\overline{n}).$\thinspace \ Extending the
definition (\ref{E.n}) of good historical paths, we write $\,E_{N-1}$%
\thinspace \ for the set of paths $\,w$\thinspace \ such that 
\begin{equation}
\int_0^{T_N}L_{[w,\varrho ]}({\rm d}s)\,\ge \,\xi _{N-1}\,:=\,\frac{\alpha
\,2^{N/2}}6\,.
\end{equation}
Let $\,\overline{T}_N:=2^N/6.$\thinspace \ Recall that $\,T_N=\overline{T}%
_N\wedge \tau _{N\,},$\thinspace \ and note that $\,\xi _{N-1}=\alpha \,%
\overline{T}_N\,\overline{\zeta }_{N\,}.$\thinspace \ Then 
\begin{equation}
\hspace*{-10pt}\left. 
\begin{array}{lll}
\Pi _{0,0}\left( E_{N-1}^{{\rm c}}\right)  & \le  & \Pi _{0,0}\bigg(%
\dint\nolimits_0^{T_N}\!L_{[W,\varrho ]}({\rm d}s)\,\le \,\alpha \,\overline{%
T}_N\overline{\,\zeta }_{N\,}\bigg)\vspace{8pt} \\ 
& \le  & \Pi _{0,0}\left( \tau _N\le \overline{T}_N\right) \vspace{8pt} \\ 
&  & +\;\Pi _{0,0}\bigg(\dint\nolimits_0^{T_N}\!L_{[W,\varrho ]}({\rm d}%
s)\,\le \,\alpha \,\overline{T}_N\,\overline{\zeta }_{N\,},\,\;\tau _N>%
\overline{T}_N\bigg)\vspace{8pt} \\ 
& \le  & \Pi _{0,0}\left( \tau _N\le \overline{T}_N\right) \vspace{8pt} \\ 
&  & +\;\Pi _{0,0}\bigg(\dint\nolimits_0^{\overline{T}_N}\!L_{[W,\varrho ]}(%
{\rm d}s)\,\le \,\alpha \,\overline{T}_N\,\overline{\zeta }_N\bigg).
\end{array}
\right\}   \label{previous.estimate}
\end{equation}
Let $\,\overline{\varrho }=\overline{\varrho }(N)$\thinspace \ be obtained
from $\,\varrho $\thinspace \ as follows. Let $\,\overline{\varrho }%
_b:=\varrho _b$\thinspace \ if $\,b\in D_{N\,}.$\thinspace \ Otherwise, set $%
\,\overline{\varrho }$\thinspace $_b:=1.$\thinspace \ \thinspace Let $\,S_N$%
\thinspace \ denote the collection of sites $\,b\in {\sf Z}^d$\thinspace \
such that $\,\overline{\varrho }\le \overline{\,\zeta }_{N\,}$\thinspace \
[recall notation (\ref{notation.zeta.bar})]. Note that by the definition of $%
\,A\left( m,N,\overline{\zeta }_N\right) \supseteq A_2(\overline{n}),$%
\thinspace \ there is no connected subset of $\,S_N$\thinspace \ with
cardinality $\,N.$\thinspace \ Thus, by Lemma \ref{lemma.clusters}, 
\[
\Pi _{0,0}\left( \int_0^{\overline{T}_N}\!{\rm d}s\;{\bf 1}\left\{ W_s\in
S_N^{{\rm c}}\right\} \;\le \;\alpha \overline{T}_N\right) \,\le
\;c^{-1}\,\exp \left[ -c\,\overline{T}_N\right] .
\]
By the definition of $\,S_{N\,}$, 
\[
\int_0^{\overline{T}_N}\!{\rm d}s\;{\bf 1}\left\{ W_s\in S_N^{{\rm c}%
}\right\} \;>\;\alpha \overline{T}_N\quad {\rm implies}\quad \int_0^{%
\overline{T}_N}\!L_{[W,\overline{\varrho }]}({\rm d}s)\;>\;\alpha \overline{T%
}_N\,\overline{\zeta }_{N\,},
\]
and so 
\[
\begin{array}{l}
\Pi _{0,0}\left( \dint\nolimits_0^{\overline{T}_N}\!L_{[W,\overline{\varrho }%
]}({\rm d}s)\;\le \;\alpha \overline{T}_N\,\overline{\zeta }_N\right) 
\vspace{8pt} \\ 
\le \;\Pi _{0,0}\left( \dint\nolimits_0^{\overline{T}_N}\!{\rm d}s\;{\bf 1}%
\left\{ W_s\in S_N^{{\rm c}}\right\} \;\le \;\alpha \overline{T}_N\right) 
\end{array}
\]
Therefore, 
\[
\Pi _{0,0}\left( \int_0^{\overline{T}_N}\!L_{[W,\overline{\varrho }]}({\rm d}%
s)\;\le \;\alpha \overline{T}_N\,\overline{\zeta }_N\right) \;\le
\;c^{-1}\exp \left[ -c\,\overline{T}_N\right] .
\]
Now (\ref{previous.estimate}) and the previous estimate combined with (\ref
{3.14}) gives 
\begin{eqnarray*}
\Pi _{0,0}\left( E_{N-1}^{{\rm c}}\right)  &\le &\Pi _{0,0}\left( \tau _N\le 
\overline{T}_N\right) +c^{-1}\exp \left[ -c\,\overline{T}_N\right]  \\
&\le &\left( \frac{2^N\,\frac{{\rm e}}6}{2^N}\right) ^{\!2^N}\left( 2\pi
2^N\right) ^{-\frac 12}+c^{-1}\exp \left[ -c\,\,\frac{2^N}6\right]  \\%
[0.05cm]
&\le &\exp \left[ -2c_3\,2^N\right] 
\end{eqnarray*}
(increasing $\,N$\thinspace \ if necessary). Therefore, by Markov's
inequality, and the `stopped expectation' formula (\ref{stopped.expectation}%
), 
\begin{equation}
\widetilde{P}_{0,\delta _0}^\varrho \left( \widetilde{X}_{T_N}^\varrho
(E_{N-1}^{{\rm c}})>\exp \left[ -c_3\,2^N\right] \right) \,\le \;\exp \left[
-c_3\,2^N\right] .  \label{Carl}
\end{equation}

Next we consider $\,E_{N-1\,}.$\thinspace \ We wish to show that in the case 
$\,K=L_{[W,\varrho ]\,},$\thinspace \ condition (\ref{infinite.BCLT}) in
Proposition \ref{P.comparison.Feller}\thinspace \ holds $%
\bP%
$--a.s. By Fubini's theorem, it suffices to verify it $\,\Pi _{0,0}\times 
\bP%
$--a.s. First note that with $\Pi _{0,0}$--probability 1 the range ${\cal R}%
(W)$\thinspace \ of the random walk $\,W$\thinspace \ is infinite. For each
site $\,b\in {\cal R}(W),$\thinspace \ let $\,Y_b\,:=\,\sigma _{b\,}\varrho
_{b\,},$\thinspace \ where $\,\sigma _b$\thinspace \ is the amount of time
which $\,W$\thinspace \ spends at $\,b$\thinspace \ between the time of
first arrival at $\,b$\thinspace \ and the first subsequent departure. Then
the $\,Y_b$\thinspace \ are i.i.d.\ with positive $\Pi _{0,0}\times 
\bP%
$--expectation. Therefore, by the strong law, 
\[
\int_0^\infty \!L_{[W,\varrho ]}({\rm d}s)\;\ge \;\sum_{b\in {\cal R}%
(W)}Y_b\;=\;\infty ,\quad \Pi _{0,0}\times 
\bP%
{\rm -a.s.,}
\]
giving (\ref{infinite.BCLT}).

By Proposition \ref{P.comparison.Feller} with 
\[
T_0=0,\quad T_1=T_{N-1\,},\quad \xi =\xi _{N-1}=\left( \alpha /6\right)
\,2^{N/2},\quad {\rm and}\quad E=E_{N-1\,}, 
\]
we obtain 
\begin{equation}
\widetilde{P}_{0,\delta _0}^\varrho \left( \widetilde{X}_{T_{N_{\!_{}}}}^%
\varrho (E_{N-1})>0\right) \,\le \,\frac 1{\xi _{N-1}}\,=\,\frac 6\alpha
\,2^{-N/2}.  \label{Carl.2}
\end{equation}

If $\,N\,\,\;$is large enough, we have $\,M_N\ge \exp \left[
-c_3\,2^N\right] .$\thinspace \ Hence, by (\ref{Carl}) and (\ref{Carl.2}), 
\begin{eqnarray*}
\widetilde{P}_{0,\delta _0}^\varrho \left( \big\| \widetilde{X}%
_{T_N}^\varrho \big\| >M_N\right)  &\le &\widetilde{P}_{0,\delta _0}^\varrho
\left( \widetilde{X}_{T_N}^\varrho \left( E_{N-1}^{{\rm c}}\right) >\exp
\left[ -c_3\,2^N\right] \right)  \\[2pt]
&&+\;\widetilde{P}_{0,\delta _0}^\varrho \left( \widetilde{X}_{T_N}^\varrho
\left( E_{N-1}\right) >0\right)  \\
&\le &\exp \left[ -c_3\,2^N\right] +\frac 6\alpha \,2^{-N/2}\;\,\le
\,\;2^{-N/4}\;\,=\,\;\delta _{N-1\,}.
\end{eqnarray*}
This implies the starting condition (\ref{starting.cond}) in Hypothesis \ref
{H.good.and.bad}.

Now we consider the other time intervals $\,[T_{n\,},T_{n+1}),$\thinspace $%
\;n\ge N.$\thinspace \ To deal with these times, we no longer consider
clusters of sites where the catalyst is small, but just consider single
sites. Recall that we are on the set $\,A_1(\overline{n})$\thinspace \ and
that $\,n\ge N\ge \overline{n}.$\thinspace \ Note that on $%
\,A_1(n+1)\supseteq A_1(\overline{n})$\thinspace \ we have $\,\varrho
_b>\zeta _{n+1\,},$\thinspace \ for $\,b\in D_{n+1\,}.$\thinspace \
Therefore, if $\,W_s\in D_{n+1}$\thinspace \ for $\,T_n\le s\le T_n+2^{-n},$
then $\,T_{n+1}=T_n+2^{-n}$\thinspace \ and $\,$%
\[
\int_{T_n}^{T_{n+1}}\!L_{[W,\varrho ]}({\rm d}s)\;\ge \;2^{-n}\zeta
_{n+1}\;=\;\xi _n
\]
[recall (\ref{notation.zeta})]. Hence, 
\[
\Pi _{T_{n\,},a}\left( E_n^{{\rm c}}\right) \;\le \;\Pi _{T_{n\,},a}\left(
W_s\notin D_{n+1\!_{\!_{\,_{}}}}\,\;{\rm for\;some\;\,}s\in \left[
T_{n\,},T_n+2^{-n}\right] \right) ,\quad \,\;a\in D_{n\,}.
\]
The strong Markov property applied to $\,T_n$\thinspace \ gives 
\[
\Pi _{T_{n\,},a}\left( E_n^{{\rm c}}\right) \;\le \;\Pi _{0,0}\bigg(%
\sup_{s\le 2^{-n}}\,|W_s|>2^n\bigg).
\]
From our ``traveling estimate'' (\ref{3.14}), it follows that for $\,N$%
\thinspace \ large enough, $\,n\ge N$\thinspace \ implies 
\[
\Pi _{T_{n\,},a}\left( E_n^{{\rm c}}\right) \;\le \;\left( \frac{2^{-n}e}{2^n%
}\right) ^{\!2^n}(2\pi 2^n)^{-1/2}\;\,\le \,\;c\,2^{-2\cdot 2^n}\;\,\le
\,\;\lambda _n\,\frac{M_{n+1}}{M_n}\,.
\]
Thus, Lemma \ref{L.new} gives the conditional expectation estimate (\ref
{small.expectation}).

Again, for $\,0\le s\le 2^{-n},$\thinspace \ Proposition \ref
{P.comparison.Feller} yields 
\begin{equation}
\widetilde{P}_{0,\delta _0}^\varrho \left\{ \widetilde{X}_{T_{n_{\!}}}^%
\varrho (E_n)>0\;\Big| \;\big\| \widetilde{X}_{T_n^\varepsilon }^\varrho %
\big\| \,\le \,M_n\right\} \,\le \,\frac{M_n}{\xi _n}\,=\,\delta _{n\,}.
\label{3.17}
\end{equation}
This proves the good paths estimate (\ref{Feller.survival}).

So Hypothesis \ref{H.good.and.bad} is satisfied, and {\em finite time
extinction for the lattice model}\/ follows from the abstract Theorem \ref
{T.abstract.criterion}.\hfill 
\endproof%


\vfill

\hfill {\scriptsize \quad printed \today
}\newpage 


\end{document}